\documentclass{amsart}
\usepackage{graphicx} 
\usepackage[margin=1in]{geometry}
\usepackage[initials,alphabetic]{amsrefs}
\usepackage{tikz-cd}
\usepackage{ulem}

\usepackage{amsmath}
\usepackage{amsfonts}
\usepackage{amssymb}

\newcommand{\mc}{\mathcal}
\newcommand{\mb}{\mathbf}
\newcommand{\mbb}{\mathbb}

\renewcommand{\bar}[1]{\overline{#1}}
\newcommand{\dual}{\vee}
\newcommand{\floor}[1]{\left\lfloor#1\right\rfloor}
\newcommand{\ceil}[1]{\left\lceil#1\right\rceil}

\newcommand{\ZZ}{\mathbb{Z}}

\newcommand{\CC}{\mathbb{C}}

\newcommand{\PP}{\mathbb{P}}
\newcommand{\cO}{\mathcal{O}}
\newcommand{\Pic}{\operatorname{Pic}}
\newcommand{\Sym}{\mathrm{Sym}}
\newcommand{\Hom}{\operatorname{Hom}}

\newcommand{\cF}{\mathcal{F}}

\newcommand{\boxpow}[2]{#1{}^{\boxtimes #2}}
\newcommand{\bLambda}{\mathbf{\Lambda}}
\newcommand{\cD}{\mathcal{D}}
\newcommand{\cP}{\mathcal{P}}

\newcommand{\I}{\mb{I}}
\newcommand{\II}{\mb{II}}
\newcommand{\ii}{\mb{ii}}
\newcommand{\III}{\mb{III}}
\newcommand{\iii}{\mb{iii}}
\newcommand{\IV}{\mb{IV}}
\newcommand{\iv}{\mb{iv}}
\usepackage{amsthm}
\usepackage{enumitem}
\usepackage{hyperref}
\usepackage[capitalize]{cleveref}

\numberwithin{equation}{section}

\newtheorem{thm}{Theorem}[section]
\newtheorem{prp}[thm]{Proposition}
\newtheorem{lem}[thm]{Lemma}
\newtheorem{cor}[thm]{Corollary}
\newtheorem{cnj}[thm]{Conjecture}

\theoremstyle{definition}

\newtheorem{rem}[thm]{Remark}
\newtheorem{ntn}[thm]{Notation}

\theoremstyle{remark}
\newtheorem{clm}[thm]{Claim}

\newlist{thmlist}{enumerate}{1}
\newlist{lemlist}{enumerate}{1}
\newlist{prplist}{enumerate}{1}
\newlist{corlist}{enumerate}{1}
\newlist{cnjlist}{enumerate}{1}
\newlist{dfnlist}{enumerate}{1}
\newlist{remlist}{enumerate}{1}
\newlist{ntnlist}{enumerate}{1}
\setlist[thmlist,lemlist,prplist,corlist,cnjlist,dfnlist,remlist,ntnlist]{label=\textup{(\alph*)},ref=\thethm(\alph*)}

\crefname{section}{Section}{Sections}
\crefname{subsection}{Subsection}{Subsections}

\crefname{figure}{Figure}{Figures}

\crefname{thm}{Theorem}{Theorems}
\crefname{thmlisti}{Theorem}{Theorems}

\crefname{lem}{Lemma}{Lemmas}
\crefname{lemlisti}{Lemma}{Lemmas}

\crefname{prp}{Proposition}{Propositions}
\crefname{prplisti}{Proposition}{Propositions}

\crefname{cor}{Corollary}{Corollaries}
\crefname{corlisti}{Corollary}{Corollaries}

\crefname{cnj}{Conjecture}{Conjectures}
\crefname{cnj}{Conjecture}{Conjectures}

\crefname{dfn}{Definition}{Definitions}
\crefname{dfnlisti}{Definition}{Definitions}

\crefname{ntn}{Notation}{Notations}
\crefname{ntnlisti}{Notation}{Notations}

\crefname{rem}{Remark}{Remark}
\crefname{remlisti}{Remark}{Remark}

\crefname{clm}{Claim}{Claims}

\title{Noncommutative resolution of $\mc SU_C(2)$}
\author{Elias Sink and Jenia Tevelev}
\date{\today}

\begin{document}

\begin{abstract}
    We study the derived category of the moduli space $\mc SU_C(2)$ of rank $2$ vector bundles on a smooth projective curve $C$ of genus $g\ge 2$ with trivial determinant. This generalizes the recent work by Tevelev and Torres on the case with fixed odd determinant. Since $\mc SU_C(2)$ is singular, we work with its noncommutative resolution of singularities constructed by P\u adurariu and \v Spenko--Van den Bergh (in the more general setting of symmetric stacks). We show that this noncommutative resolution admits a semiorthogonal decomposition into derived categories of  symmetric powers $\Sym^{2k}C$ for $2k\le g-1$. In the case of even genus, each block appears four times. This is also true in the case of odd genus, except that the top symmetric power $\Sym^{g-1}C$ appears twice. 
    In the case of even genus, the noncommutative resolution is strongly crepant in the sense of Kuznetsov and categorifies the intersection cohomology of $\mc SU_C(2)$.
    Since all of its components are ``geometric,'' our semiorthogonal decomposition provides evidence for the expectation, which dates back to the work of Newstead and Tyurin, that $\mc SU_C(2)$ is a rational variety.
    Finally, we study mutations of semiorthogonal decompositions on the Hecke correspondence, answering a question of 
    P\u adurariu and Toda.
\end{abstract}

\maketitle
\section{Introduction}

Let $C$ be a smooth complex projective curve of genus $g\ge2$ 
and let 
$\mc SU_C(2)$ be the coarse moduli space of semistable rank $2$ vector bundles on $C$ with fixed determinant $\Lambda$ of even degree. By tensoring with a fixed line bundle, it is easy to see that  $\mc SU_C(2)$ is independent (up to isomorphism) of the choice of $\Lambda$. Common choices are 
$\mc O_C$ or $\omega_C$, but it will be more convenient for us to choose an arbitrary $\Lambda$ such that $\deg\Lambda=2g$. It is well-known that $\mc SU_C(2)$ is a Gorenstein Fano variety of dimension $3g-3$ with rational singularities \cite{drezetnarasimhanpicard}. Various resolutions of singularitities of $\mc SU_C(2)$ have been studied in \cites{seshadridesing,narasimhanramanangeometry,kirwanhomology}, and the relationships between these desingularizations have been worked out in \cites{kiemlidesing,choechoykiemhecke}.
On the other hand, Kuznetsov \cite{kuznetsovcatres}
defines a \textit{noncommutative resolution of singularities} of a variety $X$ as a smooth triangulated category $\mc D$ with an adjoint pair of functors $f_*:\,{\mc D}\to D^b(X)$ and $f^*:\,\mathrm{Perf}(X)\to {\mc D}$
such that $f_*\circ f^*\cong\mathrm{Id}$. As our varieties are all projective, we also require $\mc D$ to be proper. In \cite{kuznetsovluntscatres}*{Theorem 1.4}, it is shown that every variety $X$ admits such a noncommutative resolution, which is proper if $X$ is.

The main result of this paper is the following theorem, which provides an even-degree counterpart 
of the main result in \cites{tevelevtorresbgmn,tevelevbraid}
(the proof of the BGMN conjecture).

\begin{thm}\label{thm:ncr}
    There exists a noncommutative resolution of singularities $\mc D$ of $\mc SU_C(2)$ with a semiorthogonal decomposition into blocks equivalent to $D^b(\Sym^{2k}C)$ for $2k\le g-1$. There are four copies of each block except when $g$ is odd, in which case the block  $D^b(\Sym^{g-1}C)$ appears twice. The category $\mc D$ is an example of the noncommutative resolution of singularities constructed in \cites{padurariu,spenko} for symmetric stacks. In particular, it is  an admissible subcategory of the derived category of the Kirwan resolution of $\mc SU_C(2)$.
 If~$g$~is~even,  $\mc D$ is a strongly crepant noncommutative resolution of  $\mc SU_C(2)$ in the sense of \cite{kuznetsovcatres}.
\end{thm}

We refer to Theorem~\ref{thm:main_precise} 
for a more concise and detailed statement of the main theorem.
As in \cites{tevelevtorresbgmn,tevelevbraid}, we construct $\mc D$ as an admissible subcategory in $D^b(M)$, where $M$ is the Thaddeus moduli space of stable pairs $(F,s)$ with $F$ a rank $2$ vector bundle on $C$ with determinant $\Lambda$ and $s$ a non-zero global section  (see \cite{thaddeusstablepairs}). Furthermore, derived categories $D^b(\Sym^{2k}C)$ are embedded in $D^b(M)$ by means of explicit Fourier--Mukai functors with kernels given by ``tensor'' vector bundles twisted by various line  bundles. 

An easy application of the main result in \cite{padurariu} (see Remark~\ref{rem:motive}) shows that $\mc D$ categorifies
the intersection cohomology of $\mc SU_C(2)$
in the even genus case. Indeed, our decomposition is compatible with its computation in \cite{delbanomotive}.
By \cite{kuznetsovcatres}*{Conjecture 4.10}, we expect $\mc D$ to be minimal (i.e., admissible in every noncommutative resolution of $\mc SU_C(2)$) when $g$ is even. Bondal and Orlov \cite{BOderived} have conjectured that every variety admits a minimal noncommutative resolution; see also \cite{halpernleistnerncmmp}.

In the odd degree case ($\deg \Lambda=2g-1$) studied in \cites{tevelevtorresbgmn,tevelevbraid}, the moduli spaces of stable bundles and stable pairs are birational.
In contrast, when $\deg \Lambda=2g$, the morphism $M\to{\mc SU_C(2)}$ has generic fibers~$\mbb P^1$. 
While $M$ is a rational variety, the rationality of the moduli space 
$\mc SU_C(2)$ remains a well-known open problem, dating back to the early works of 
Tyurin \cites{tyurin1,tyurin2}
and Newstead \cites{newstead, newsteadcorr}. Unlike the case of coprime rank and degree---where rationality has been completely settled in \cite{kingschofield}---rationality for 
$\mc SU_C(2)$ is only known in genus $2$, where the theta-morphism \cite{raynaud}
$\mc SU_C(2)\to \mathbb P^{2^g-1}$ happens to be an isomorphism~\cite{narasimhanramananriemannsurface}. In light of the ``Kuznetsov rationality proposal'' \cite{kuznetsov2016}, our main result suggests that any weak factorization of a hypothetical birational map $\mc SU_C(2)\dashrightarrow{\PP}^{3g-3}$ should involve blow-ups and blow-downs of even symmetric powers $\Sym^{2k}C$ for $2k\le g-1$.

In the even genus case, \cref{thm:ncr} proves a conjecture of Belmans \cite{belmans}. Although in odd genus it differs from that prediction, the latter can be salvaged by considering twisted noncommutative resolutions, namely the \textit{quasi-BPS categories} $\mbb{B}_w$ introduced by P\u adurariu and Toda in \cite{padurariutodaquasibps}*{Section 3.4}.
These are subcategories of the derived category of the stack
of semistable rank $2$ vector bundles with trivial determinant. In particular, $\mbb{B}_0$ is the noncommutative resolution $\cD$ from \cref{thm:ncr}. In \cref{thm:hecke_sod_cd}, we construct an analogous semiorthogonal decomposition of $\mbb{B}_1$ into blocks equivalent to $D^b(\Sym^{2k+1}C)$ for $2k+1\le g-1$.

There is  a large body of conjectures on explicit semiorthogonal decompositions of Fano varieties, which are often obtained by analyzing  Hodge diamonds and using the additivity of Hochschild homology in semiorthogonal decompositions  \cite{keller}*{Theorem~1.5}; see \cite{belmans2024} for a recent example related to our paper.\break 
These predictions are not easy to turn into proofs. 
One general framework 
is to analyze a ``two-ray game" given by two extremal contractions of Fano varieties 
\begin{equation}\label{eq:tworay}
\begin{tikzcd}[column sep=small]
& X \arrow[dl, dashrightarrow] \arrow[dr, dashrightarrow] & \\
Y  & & Z
\end{tikzcd}
\end{equation}
In our paper, $Y=\mathbb P^{3g-2}$, $X=M$, and $Z=\mc SU_C(2)$. We start with a known semiorthogonal decomposition of $D^b(Y)$ and ``extend'' it to a semiorthogonal decomposition of $D^b(X)$.
The difficult part is to mutate this decomposition of $D^b(X)$ into a decomposition compatible with a semi-orthogonal decomposition of $D^b(Z)$:


\begin{cnj}[Two-Ray Game Conjecture (1)] \label{cnj:tworay}
Given extremal contractions \eqref{eq:tworay}
of smooth Fano varieties, there exist semiorthogonal decompositions 
``compatible'' with pullbacks along maps $Y\dashleftarrow X\dashrightarrow Z$:
$$D^b(Y)=\langle\mc A_1,\ldots,\mc A_s\rangle,
\ 
D^b(X)=\langle \mc A_1,\ldots,\mc A_s, \mc P_1,\ldots,\mc P_r\rangle
=\langle Q_1,\ldots,\mc Q_u, \mc B_1,\ldots,\mc B_t\rangle,\ 
D^b(Z)=\langle\mc B_1,\ldots,B_t\rangle.$$
The two semiorthogonal decompositions
of $D^b(X)$ are related by a mutation. 
\end{cnj}

Our~approach uses weaving patterns of \cite{tevelevbraid}, allowing for tight control of the Fourier--Mukai kernels for the various functors
$D^b(\Sym^kC)\to D^b(M)$
embedding the blocks on different stages of the mutation. 
From~the~perspective of homological mirror symmetry for Fano manifolds
(see the survey \cite{auroux}),
such mutations should exist in general, and we can even  
predict a description of the hypothetical braid:

\begin{cnj}[Two-Ray Game Conjecture (2)]
The  braid in Conjecture~\ref{cnj:tworay} can be computed as the monodromy braid of eigenvalues of $c_1(X)$ acting on quantum cohomology $\mathbb QH^*(X,\mathbb C)$ as the base $\tau$ of small quantum cohomology varies along the path in the ample cone of $X$ (with a small $B$-field perturbation $iB$ added to the path $\tau$ to avoid collisions of eigenvalues). When the path approaches the walls of the ample cone, eigenvalues agglomerate in groups that depend on the structure of the birational contractions $Y\dashleftarrow X\dashrightarrow Z$.
\end{cnj}

At the moment, the only way to verify this conjecture is to compute and compare the two braids; see~
\cite{iritani}*{Figure 16} for a worked-out example
of extremal contractions $\mathbb P^4\leftarrow\text{Bl}_{\mathbb P^1}\mathbb P^4\rightarrow\mathbb P^2$ or
\cite{genus2video} for contractions  $\mathbb P^3\leftarrow \text{Bl}_{C}\mathbb P^3=\text{Bl}_{\mathbb P^1} Q_4\rightarrow Q_4$,
where $C\hookrightarrow \mathbb P^3$
is a quintic curve of genus $2$ and $B_4$ is the intersection of two quadrics in $\mathbb P^5$
(this is the genus $2$ case of 
\cite{tevelevbraid}.) 

It seems plausible that smoothness assumptions in Conjecture~\ref{cnj:tworay}
can be weakened by considering noncommutative resolutions of singularities 
as in our paper: the moduli space
$\mc SU_C(2)$ is singular,
but most of the blocks in the semiorthogonal decomposition of $D^b(M)$ ``fly away" to leave only a noncommutative resolution of  $\mc SU_C(2)$ (see Figure~\ref{fig:genus5plainweave}).
It would be interesting to gather more evidence supporting this conjecture by connecting other Fano varieties $Y$ and $Z$ through a Fano variety $X$. One could consider toric Fano varieties, maximal flag varieties, Fano $3$-folds, and various Fano moduli spaces.

Following a suggestion of P\u adurariu and Toda, we  study Conjecture~\ref{cnj:tworay} for the Hecke correpondence
$$
\begin{tikzcd}[column sep=small]
& P \arrow[dl] \arrow[dr, "\pi"] & \\
\mc SU_C(2)  & & \hat N
\end{tikzcd}
$$
Here $\hat N$ is the moduli space of stable rank $2$ vector bundles on $C$ with fixed odd determinant
and $P$ is the moduli space of stable parabolic bundles (at a fixed point $q\in C$).
We have semiorthogonal decompositions
$$\langle\mbb B_0,\mbb B_1\rangle=D^b(P)=\langle D^b(\hat N), D^b(\hat N)\otimes\mc O_\pi(1)\rangle.$$
Here $\mbb B_0=\cD$ and $\mbb B_1$ are the quasi-BPS categories, while the second 
semiorthogonal decomposition reflects the fact that $\pi$ is a $\mbb P^1$-bundle.
According to \cites{tevelevbraid, tevelevtorresbgmn}, the right-hand side admits a semiorthogonal decomposition into derived categories of symmetric powers $\Sym^k C$, as does the left-hand side by the results of this paper. In \cref{thm:hecke_sod_cd}, we construct a mutation relating these two semiorthogonal decompositions, which we call the Hecke Braid.

The paper is organized as follows: In \cref{sec:stablepairssod}, we obtain several explicit semiorthogonal decompositions of the ``maximal" (see \cref{rem:maximal}) Thaddeus space $M_{i_d}(d)$ for various $d$ (see, e.g., \cref{thm:sod} for a precise statement), including the case $d=2g$ we need (\cref{thm:sod_2g}). This is done by generalizing the weaving techniques developed in \cite{tevelevbraid} for $d=2g-1$. This in turn requires the careful verification of a number of technical results from \cite{tevelevbraid} for other degrees; these checks are carried out in \cref{sec:basicweaving}. The precise statement and proof of the main result are given in \cref{sec:plainweave}. In \cref{sec:hecke}, we study the Hecke Braid.

A few words regarding notation: Following \cites{thaddeusstablepairs,tevelevbraid}, we often denote tensor product by juxtaposition for compactness. As in \cite{huybrechtsfouriermukai}, we usually omit $R$'s and $L$'s on derived functors except for emphasis (e.g., when applying derived pushfoward to a sheaf). We frequently use the same symbol to denote canonical objects on related moduli spaces when no confusion will arise, omitting explicit pullbacks (see \cref{ntn:stablepairs}).

\subsection*{Acknowledgements}
We are grateful to Igor Dolgachev for the suggestion to study the derived category of $\mc SU_C(2)$ in relation to its conjectural rationality, to James Hotchkiss for help with  derived categories of stacks, to Sasha Kuznetsov and Tudor P\u adurariu for insightful discussions, and to participants of the Spring 2024 seminar on the noncommutative minimal model program at UMass Amherst for  an inspirational research environment. 
The second author would like to thank Yukinobu Toda for discussions and kind hospitality during a visit to Kavli IPMU.
This research was  supported by NSF grants DMS-2101726 and DMS-2401387.
\section{Semiorthogonal decompositions of moduli spaces of stable pairs}\label{sec:stablepairssod}

We begin by recalling some notation from \cites{tevelevbraid,tevelevtorresbgmn}.
\begin{ntn}\label{ntn:stablepairs}
    For a line bundle $\Lambda$ on $C$ of degree $d$ and $0\leq i\leq \floor{\frac{d-1}{2}}$, we denote by $M_i(\Lambda)$ (or simply $M_i(d)$ or $M_i$ when no confusion will arise) the moduli space of rank 2 stable pairs with determinant $\Lambda$, where $i$ indexes the stability condition. We write $\cO(m,n)=\cO((m+n)H-nE)$ on any of the $M_i$, $i\geq 1$, where $E$ is the exceptional divisor of the contraction $M_1\to M_0$ and $H$ is the pullback of the hyperplane divisor from $M_0\cong \PP^{d+g-2}$. (By abuse of notation, $\cO(m,0)=\cO(m)$ on $M_0$.) $\mc F$ denotes the universal vector bundle on $M_i\times C$ (for any $i$), $\mc F_x$ its restriction to $M_i\times \{x\}\cong M_i$ for $x\in C$, and $\bLambda$ the line bundle $\wedge^2\mc F_x$, which is independent of $x$ (and not to be confused with $\Lambda$). 
    
    For any variety $X$ and vector bundle $\mc G$ on $X\times C$, we have tensor vector bundles $\boxpow{\mc G}{k}$ and $\boxpow{\bar{\mc G}}{k}$ on $X\times \Sym^k C$ defined by the $S_k$-equivariant pushforwards $\tau_*^{S_k}(\pi_1^*\mc G\otimes\cdots \otimes \pi_k^*\mc G)$ and $\tau_*^{S_k}(\pi_1^*\mc G\otimes\cdots \otimes \pi_k^*\mc G\otimes \mathrm{sgn})$, respectively, where $\tau:C^k\to \Sym^k C$ is the quotient by $S_k$, $\pi_i:C^k\to C$ are the projections, $\mathrm{sgn}$ is the sign character of $S_k$, and the $S_k$-action on $\pi_1^*\mc G\otimes\cdots \otimes \pi_k^*\mc G$ permutes the tensor factors (see \cite{tevelevtorresbgmn}*{Section~2}). For $D\in \Sym^k C$, we denote by $\boxpow{\mc G}{k}_D$ and $\boxpow{\bar{\mc G}}{k}_D$ the restrictions to $X\times\{D\}\cong X$ of $\boxpow{\mc G}{k}$ and $\boxpow{\bar{\mc G}}{k}_D$, respectively. $\langle K\rangle$ denotes the essential image of a fully faithful Fourier--Mukai functor with kernel $K$.
\end{ntn}

The principal goal of this section is to prove the following generalization of \cite{tevelevbraid}*{Theorem 3.1}:
\begin{thm} \label{thm:sod}
    Let $d\leq 2g$ with $i_d\leq \floor{\frac{d-1}{2}}$, where $i_d=\ceil{\frac{d+g-1}{3}}-1$. Let $m=d+g-1-3i_d\in\{1,2,3\}$, and let $m_n=1$ if $m\leq n$ or $0$ otherwise. Then
    \begin{equation}\label{eq:sod1}
        D^b(M_{i_d}(d))=\left\langle
        \left\langle\bLambda^{-k}\boxpow{\cF^\dual}{j} \right\rangle_{\substack{j+k\leq i_d-m_2\\j,k\geq 0}},
        \left\langle T_1\bLambda^{-k}\boxpow{\cF^\dual}{j}  \right\rangle_{\substack{j+k\leq i_d-m_1\\j,k\geq 0}},
        \left\langle T_2\bLambda^{-k}\boxpow{\cF^\dual}{j} \right\rangle_{\substack{j+k\leq i_d\\j,k\geq 0}}
        \right\rangle
    \end{equation}
    where the blocks within each ``megablock" are ordered first by decreasing $k$, then by decreasing $j$. Here, $T_1=\cO(1,i_d-m_2)$, and $T_2=\cO(2,2i_d-m_2-m_1)$.
\end{thm}

\begin{rem}
    The assumptions of the theorem are equivalent to $d=2g-\alpha$ for $\alpha\in\{0,1,2,3,5\}$, where $i_d=g-\ceil{\frac{\alpha+2}{3}}$. The restriction $d\leq 2g$ could be removed by verifying \cref{cnj:hardervanishing} below, in which case the theorem would hold in all but finitely many degrees.
\end{rem}
\begin{rem}\label{rem:maximal}
    It can be seen from \cite{tevelevtorresbgmn}*{Proposition 3.18} that $D^b(M_i(d))$ is ``largest" when $i=i_d$. Moreover, it follows from \cite{thaddeusstablepairs}*{5.3, 6.1} that $M_{i_d}$ is Fano when $m=1,2$. When $m=3$, $D^b(M_{i_d})$ and $D^b(M_{i_d+1})$ are equivalent, and the anticanonical bundles on $M_{i_d}$ and $M_{i_d+1}$ are big and nef but not ample. 
\end{rem}

\subsection{Generalized weaving}
We begin with some notation. Fix $d\leq 2g$ and write $M_i(d)=M_i$ for $i\leq\floor{\frac{d-1}{2}}$.
\begin{ntn}\label{ntn:ft}
    For $0\leq k\leq i$, we denote by $\cD^k_i$ the structure sheaf of the reduced subscheme $$D^k_i=\{(D,F,s)\in \Sym^k C\times M_i : s|_D=0\},$$ whose fiber over $D\in \Sym^k C$ is $M_{i-k}(\Lambda(-2D))$ \cite{tevelevtorresbgmn}*{Remark 3.7}.  For $t\in[0,i_d+1)$, let $\cD^{k,s}_t=\cD^k_{\floor{t}}\otimes L^{k,s}_t$ where
    $$L^{k,s}_t=
    \begin{cases}
        \cO(s,sk)                                                          &  k = \floor{t} \\
        \cO\left(\floor{\frac{s}{t-k}},s+\floor{\frac{s}{t-k}}(k-1)\right) &  k < \floor{t}. 
    \end{cases}$$ 
\end{ntn}
Our first step is to prove the following:
\begin{lem}[cf. \cite{tevelevbraid}*{Corollary 2.10}]\label{lem:ft_main}
    For $t\in (0,i_d+1)\smallsetminus\ZZ$, we have a semiorthogonal decomposition  
    \begin{equation}\label{eq:ft_main1}
        D^b(M_{\floor{t}})=\left\langle \cD^{k,s}_t \right\rangle_{\substack{0\leq k\leq\floor{t}\\0\leq s\leq d+g-3k-2} }
    \end{equation}
    where the blocks are ordered first by increasing $x_{k,s}(t)=\frac{s}{t-k}$, then by increasing $k$.
\end{lem}
We interpret $t$ as time, with the block $\langle \cD^{k,s}_t\rangle $ ``moving" in the $x$-$t$ plane with trajectory $x=x_{k,s}(t)$ (or $x=k\epsilon$ for $s=0$, where $\epsilon\ll 1$). When the blocks cross paths, they change order and undergo mutations dictated by the line bundle $L^{k,s}_t$. When $t$ crosses an integer level $i$, we embed $D^b(M_{i-1})$ into $D^b(M_i)$, introducing several new blocks as its orthogonal complement, and proceed with the process. We refer to this ``weave" as the Farey Twill; see \cite{tevelevbraid}*{Section 2} for detailed illustrations in the case $d=2g-1$. This program is facilitated by several technical lemmas, whose statements and proofs are deferred to \cref{sec:basicweaving}. 

To pass from $M_{i-1}$ to $M_i$, we need the following ``windows" embeddings, in the sense of \cite{halpernleistnerwindows}: 
\begin{prp}[\cite{tevelevtorresbgmn}*{Proposition 3.18}]\label{prp:windows}
    For $d>0$ and $1\leq i\leq i_d\leq \floor{\frac{d-1}{2}}$, there is an admissible embedding $\iota:D^b(M_{i-1}(d))\hookrightarrow D^b(M_{i}(d))$ giving rise to a semiorthogonal decomposition
    $$D^b(M_i(d))=\langle \iota (D^b(M_{i-1}(d))),\cD^{i,0}_i,\cD^{i,1}_i,\dots,\cD^{i,d+g-3i-2}_i\rangle$$
    When $i>1$, the embedding corresponds to the inclusion of objects with weights $[0,i)\subseteq [0,d+g-1-2i)$ with respect to the wall crossing $M_{i-1}\dashrightarrow M_i$.
\end{prp}
\begin{rem}
    Note that $D^i_i$ is isomorphic to the projective bundle $\PP W^+_i\subset M_i$ via the second projection (see \cite{tevelevtorresbgmn}*{Section 6}), and that $L^{i,s}_i=\cO(s,si)$ restricts to $\cO(s)$ on the fibers $M_0(d-2i)\cong \PP^{d+g-2i-2}$ of this bundle \cite{tevelevtorresbgmn}*{Remark 3.7}.
\end{rem}

\begin{proof}[Proof of \cref{lem:ft_main}]
    When $0<t<1$, \eqref{eq:ft_main1} is the Be\u{\i}linson collection $\langle \cO,\cO(1),\dots,\cO(d+g-2)\rangle$ on $\PP^{d+g-2}$. Given \eqref{eq:ft_main1} for $t=i+\epsilon$ with $i\in\ZZ$, $\epsilon\ll 1$, we achieve \eqref{eq:ft_main1} for all $t\in(i,i+1)$ by performing the mutations encoded in the crossings of trajectories $x_{k,s}(t)$. When blocks meet at nonintegral $x$, they only change order, meaning we must show that they are mutually orthogonal. Since the intersecting blocks are already ordered by $x_{s,k}(t-\epsilon)$, or equivalently by $k$, we need only check that $\langle \cD^{k,s}_t\rangle \subset {}^\perp\langle\cD^{k',s'}_t\rangle$ for $k<k'$, $x_{k,s}(t)=x_{k',s'}(t)$. This follows from \cref{lem:ft_orthogonal} below.
    
    Crossings at $x\in \ZZ$ have the form $\langle \cD^{k,s}_t,\cD^{k+1,s-x}_t,\dots, \cD^{\floor{t},s-\floor{t}x}_t\rangle\to \langle \cD^{\floor{t},s-\floor{t}x}_{t+\epsilon},\dots,\cD^{k+1,s-x}_{t+\epsilon},\cD^{k,s}_{t+\epsilon}\rangle$. Note that for $0\leq j<\floor{t}$, $$L^{k+j,s-jx}_t=\cO(x,s-x(k-1))$$ is indepedent of $j$, while $$L^{k+j,s-jx}_{t+\epsilon}=L^{k+j,s-jx}_t(-1,1-j).$$ Moreover, $L^{\floor{t},s-\floor{t}x}_t=L^{\floor{t},s-\floor{t}x}_{t+\epsilon}=\cO(s-\floor{t}x,(s-\floor{t}x)\floor{t})$ and $\cO(x,s-x(k-1))$ both restrict to $\cO(s)$ on the fibers of the projective bundle $D^{\floor{t}}_{\floor{t}}$, so $\langle\cD^{\floor{t},s-\floor{t}x}_t\rangle= \langle\cD^{\floor{t}}_{\floor{t}}(x,s-x(k-1))\rangle$. Hence it suffices to give a mutation
    $$\langle \cD^k_i,\dots,\cD^{i-1}_i,\cD^i_i\rangle\to \langle\cD^{i}_i,\cD^{i-1}_i(-1,2-i),\dots,\cD^k_i(-1,1-k)\rangle,$$
    as in \cref{lem:ft_mutation}.

    It remains to explain how to get from $t=i-\epsilon$ to $i+\epsilon$. By \cref{lem:ft_windows} below, we have $\iota\langle\cD^{k,s}_{i-\epsilon}\rangle=\langle\cD^{k,s}_{i}\rangle$. Hence, to go from the semiorthogonal decomposition of \cref{prp:windows} to \eqref{eq:ft_main1} with $t=i+\epsilon$, we need only move the block $\langle\cD^{i,0}_i\rangle$ into position, $i$-th from the left (we imagine this block coming horizontally from the right along $t=i$, stopping at $x=i\epsilon$). It moves past blocks with $x_{k,s}(i)\not\in\ZZ$ without changing them by \cref{lem:ft_orthogonal}, while the others undergo the mutation of \cref{lem:ft_mutation}.
\end{proof}

It turns out that \cref{lem:ft_main} is not quite the decomposition we need to proceed.
\begin{lem}\label{lem:ft_modified}
    Let $m=d+g-3i_d-1\in\{1,2,3\}$, and let $m_n=1$ if $m\leq n$ and $0$ otherwise. Then
    \begin{equation}\label{eq:ft_modified1}
        D^b(M_{i_d})=\left\langle
        \left\langle\bLambda^{-j}\cD^k \right\rangle_{\substack{j+k\leq i_d-m_2\\j,k\geq 0}},
        \left\langle T_1\bLambda^{-j}\cD^k \right\rangle_{\substack{j+k\leq i_d-m_1\\j,k\geq 0}},
        \left\langle T_2\bLambda^{-j}\cD^k\right\rangle_{\substack{j+k\leq i_d\\j,k\geq 0}}
        \right\rangle
    \end{equation}
    where the blocks in each megablock are ordered first by increasing $j+k$, then by increasing $j$. Here, $\cD^k=\cD^k_{i_d}$, $T_1=\cO(1,i_d-m_2)$, and $T_2=\cO(2,2i_d-m_2-m_1)$.
\end{lem}

\begin{proof}
    For $m=1,2$, we begin with \eqref{eq:ft_main1} with $t=i_d-\epsilon$. We embed with $\iota$ to obtain the following semiorthogonal decomposition: 
    \begin{equation}\label{eq:ft_roman}
        D^b(M_{i_d})=\langle\I,\II,\III,\IV,\cD^{i_d,0}_{i_d},\dots,\cD^{i_d,m-1}_{i_d}\rangle
    \end{equation}
    where $\I,\II,\III,\IV$ are the subcategories generated by $\langle\cD^{k,s}_{i_d}\rangle$ for $x_{k,s}(i_d)\in[0,1)$, $[1,2)$, $[2,3)$, and $[3,\infty)$, respectively.

    If $m=1$, the blocks are arranged as in \cite{tevelevbraid}. The blocks in $\I$ are the same as in the first megablock of \eqref{eq:ft_modified1} (with $j=s$), but they are ordered differently. We reorder them by moving $\langle\cD^{k,s}_{i_d}\rangle$ from $x=x_{k,s}(i_d)$ to $x=\frac{s+k}{i_d}$. The moves are done in order of decreasing $s+k$, then by decreasing $s$. As blocks in $\I$ with the same $s+k$ are already ordered by increasing $s$, the orthogonality we need to ensure no mutations occur is $\langle\cD^{k,s}_{i_d}\rangle \subset{}^\perp\langle\cD^{k',s'}_{i_d}\rangle$ for $s'+k'<s+k$. This is checked in \cref{lem:ft_reordering}. 
    

    Similarly, the blocks in $\II$ and $\III$ are respectively the same as the second and third (with $j+k\leq i_d-1$) megablocks of \eqref{eq:ft_modified1}. (Explicitly, we have $j= s-\floor{\frac{s}{i_d-k}}(i_d-k)$.) The same reordering procedure works, so we are left to produce the blocks in \eqref{eq:ft_modified1} in with $j+k=i_d$. These are exactly the blocks in $\langle \IV,\cD^{i_d,0}_{i_d}\rangle=\langle\cD^{0,3i_d}_{i_d},\cD^{1,3(i_d-1)}_{i_d},\dots,\cD^{i_d,0}_{i_d}\rangle$ after the mutation of \cref{lem:ft_mutation} (note that we have $\cD^{i_d}=\cD^{i_d}(2,2(i_d-1))$ by \cite{tevelevtorresbgmn}*{Remark 3.7}). This proves the lemma for $m=1$.

    For $m=2$, the only new blocks in \eqref{eq:ft_roman} compared to $m=1$ lie in $\IV$. We reorder the blocks in $\I$ and $\II=\II_a$ just as before; this gives the first and second (with $j+k\leq i_d-1$) megablocks in \eqref{eq:ft_modified1}. We write $\III=\langle \ii_b,\III_a\rangle$, where $\ii_b$ contains those blocks $\langle\cD^{k,s}_{i_d}\rangle$ with $x_{k,s}(i_d)=2$ and $\III_a$ those with $x_{k,s}(i_d)\in (2,3)$. Similar to the proof of \cref{lem:ft_main}, we move the block $\cD^{i_d,0}_{i_d}$ past $\IV$, $\III_a$, and $\ii_b$. By \cref{lem:ft_orthogonal,lem:ft_mutation}, $\III_a$ is unchanged, while $\IV$ and $\ii_b$ undergo some mutations. This yields
    $$D^b(M_{i_d})=\langle \I,\II_a,\II_b, \III_a,\IV',\cD^{i_d,1}_{i_d}\rangle$$
    where $\langle \ii_b, \cD^{i_d,0}_{i_d}\rangle \to \II_b$ via \cref{lem:ft_mutation}, and $\IV'=\langle\cD^{s,k}_{i_d+\epsilon}\rangle_{x_{k,s}(i_d)\geq 3}$ ordered by $x_{k,s}(i_d+\epsilon)$. Now $\II_b=\langle\cD^{k}(1,2i_d-k-1)\rangle_{0\leq k\leq i_d}$ ordered by decreasing $k$, so $\langle \II_a,\II_b\rangle$ forms the second megablock of \eqref{eq:ft_modified1}.

    At this point, $\III_a$ contains the blocks from the third megablock of \eqref{eq:ft_modified1} with $j+k\leq i_d-2$; we apply the same algorithm to put them in the correct order. We write $\IV'=\langle \III_b,\iii_c\rangle$ with $\III_b$ and $\iii_c$ containing the blocks with $x_{k,s}(i_d)=3$ and $x_{k,s}(i_d)>3$. We have $\III_b=\langle \cD^k(2,3i_d-k-2)\rangle_{0\leq k\leq i_d-1}$ ordered by decreasing $k$. The Farey Twill trajectories of blocks in $\iii_c$ meet with $\cD^{i_d,1}_{i_d}$ at $(t,x)=\left(i_d+\frac{1}{3},3\right)$, where they undergo a final mutation $\langle \iii_c,\cD^{i_d,1}_{i_d}\rangle \to \III_c=\langle\cD^k(2,3i_d-k-1)\rangle_{0\leq k\leq i_d}$ ordered by decreasing $k$. To sum up, we have
    $$D^b(M_{i_d})=\left\langle \I,\langle\II_a,\II_b\rangle,\langle\III_a,\III_b,\III_c\rangle\right\rangle,$$
    which is exactly \eqref{eq:ft_modified1} with $m=2$.

    Finally, $m=3$ is the easiest case. We begin with \eqref{eq:ft_main1} with $t=i_d+1-\epsilon$. We have $x_{k,s}(t)\in [0,3)$, where the blocks in $[0,1)$, $[1,2)$, and $[2,3)$ correspond exactly to the respective megablocks in \eqref{eq:ft_modified1} with $j=s$, $j=s-(i_d+1-k)$, and $j=s-2(i_d+1-k)$; they are in the wrong order, but this is rectified by \cref{lem:ft_reordering} and the same reordering algorithm as above.
\end{proof}

From here, \cref{thm:sod} follows exactly as in the proof of \cite{tevelevbraid}*{Theorem 3.1}.

\begin{proof}[Proof of \cref{thm:sod}]
Each megablock in \eqref{eq:ft_modified1} will mutate into the corresponding one in \eqref{eq:sod1}. As the megablocks differ only in size and overall line bundle twists (i.e., the shapes are the same), it will suffice to describe this mutation for the first one. We rely on the Cross Warp mutation depicted in \cref{fig:crosswarp} and proved as \cref{thm:crosswarp_mutation} below. Notice that the top left portion of the mutation with top center block $\cD^k$ is precisely the bottom right portion of the mutation with top center $\cD^{k-1}$; similarly, the bottom left is the top right of the mutation with top center $\cD^{k-1}$, tensored with $\bLambda^{-1}$. Hence we can stack these mutations with top centers as in \cref{fig:crosswarp_stack}. In the end, all $\cD$'s are replaced by $\cF$'s, resulting in \eqref{eq:sod1}.
\end{proof}

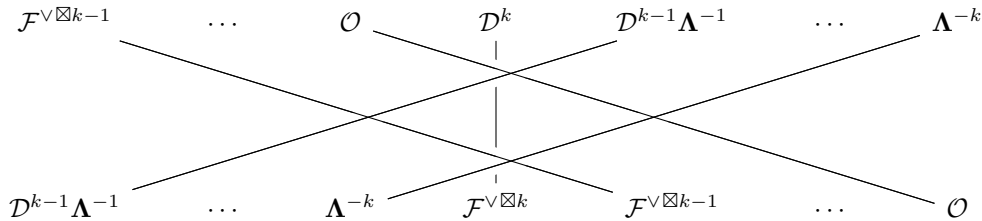
\begin{figure}[b]
    \centering
    \[\begin{tikzcd}
	{\boxpow{\cF^\dual}{k-1}} & \cdots & {\mathcal{O}} & {\mathcal{D}^k} & {\mathcal{D}^{k-1}\mathbf{\Lambda}^{-1}} & \cdots & {\mathbf{\Lambda}^{-k}} \\
	\\
	\\
	 {\mathcal{D}^{k-1}\mathbf{\Lambda}^{-1}} & \cdots & {\mathbf{\Lambda}^{-k}} & \boxpow{\cF^\dual}{k} & \boxpow{\cF^\dual}{k-1} & \cdots & {\mathcal{O}}
    \arrow[from=1-4,to=4-4,no head]
	\arrow[no head, crossing over, from=1-3, to=4-7]
	\arrow[no head, crossing over, from=1-1, to=4-5]
	\arrow[no head, crossing over, from=1-5, to=4-1]
	\arrow[no head, crossing over, from=1-7, to=4-3]
 	\arrow[no head, from=1-3, to=4-7]
	\arrow[no head, from=1-1, to=4-5]
	\arrow[no head, from=1-5, to=4-1]
	\arrow[no head, from=1-7, to=4-3]
    \end{tikzcd}\]
    \caption{The basic Cross Warp mutation, cf. \cite{tevelevbraid}*{Figure 7}.}
    \label{fig:crosswarp}
\end{figure}

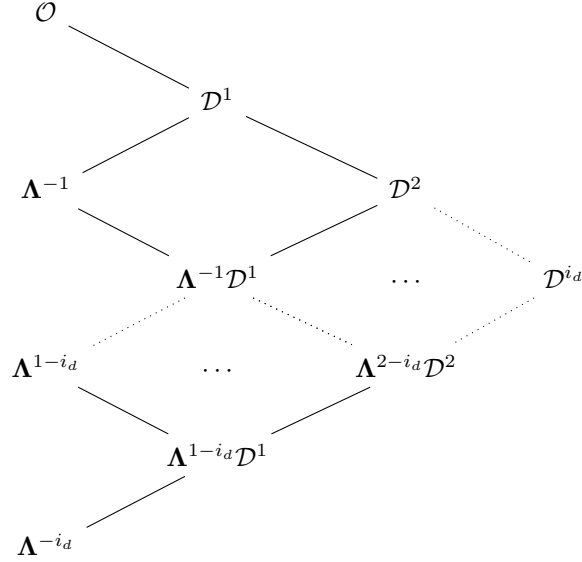
\begin{figure}
    \centering
    \[\begin{tikzcd}
	\cO \\
	& \cD^1 \\
	\bLambda^{-1} && \cD^2 \\
	& \bLambda^{-1}\cD^1 & \cdots & \cD^{i_d} \\
	\bLambda^{1-i_d} &\cdots& \bLambda^{2-i_d}\cD^2 \\
	& \bLambda^{1-i_d}\cD^1 \\
	\bLambda^{-i_d}
	\arrow[no head, from=1-1, to=2-2]
	\arrow[no head, from=2-2, to=3-1]
	\arrow[no head, from=3-1, to=4-2]
	\arrow[dotted, no head, from=4-2, to=5-1]
	\arrow[no head, from=2-2, to=3-3]
	\arrow[no head, from=3-3, to=4-2]
	\arrow[dotted, no head, from=3-3, to=4-4]
	\arrow[dotted, no head, from=4-2, to=5-3]
	\arrow[dotted, no head, from=4-4, to=5-3]
	\arrow[no head, from=5-3, to=6-2]
	\arrow[no head, from=5-1, to=6-2]
	\arrow[no head, from=6-2, to=7-1]
    \arrow[dotted, no head, from=4-2, to=5-3]
\end{tikzcd}\]
    \caption{Stacking the crosswarp mutation (with $m_1=0$). See also \cite{tevelevbraid}*{Figure 8}.}
    \label{fig:crosswarp_stack}
\end{figure}

\subsection{Broken Loom for \texorpdfstring{$d=2g$}{d=2g}}
To finish this section, we specialize to $d=2g$. We wish to modify the semiorthogonal decomposition of \cref{thm:sod} for $d=2g$ to create as many blocks as possible of the form $\theta^j\bLambda^k\boxpow{\cF^\dual}{2k}$, where $\theta=\cO(1,g-1)$ is the pullback of the ample generator of $\Pic \mc SU_C(2)$ under the forgetful morphism $M_{g-1}(2g)\to SU_C(2)$, $(F,s)\mapsto s$ (see \cite{thaddeusstablepairs}*{5.8}). We will see in \cref{sec:plainweave} that such blocks form the claimed noncommutative resolution of $SU_C(2)$.
\begin{thm}[cf. \cite{tevelevbraid}*{Theorem 5.8}]\label{thm:sod_2g}
Let $M=M_{g-1}(2g)$. We have the following semiorthogonal decomposition of $D^b(M)$:
\begin{align*}
\Bigg{\langle} 
\left\langle\theta^{-1}\bLambda^{\floor{\frac{g-2}{2}}-k}\boxpow{\cF^\dual}{\lambda-2k}\right\rangle
_{\substack{0\leq\lambda\leq g-2\\0\leq k\leq\lfloor\frac{\lambda}{2}\rfloor}},
\left\langle\bLambda^{\floor{\frac{g-2}{2}}-k}\boxpow{\cF^\dual}{\lambda-2k}\right\rangle_{\substack{0\leq\lambda\leq 2(g-2)\\0\leq k\leq\lfloor\frac{\lambda}{2}\rfloor,\lambda-k\leq g-2}},&\\
\left\langle\theta\bLambda^{\floor{\frac{g}{2}}-k}\boxpow{\cF^\dual}{\lambda-2k}\right\rangle_{\substack{0\leq\lambda\leq 2(g-1)\\0\leq k\leq\lfloor\frac{\lambda}{2}\rfloor,\lambda-k\leq g-1}},
\left\langle\theta^2\bLambda^{\floor{\frac{g}{2}}-k}\boxpow{\cF^\dual}{\lambda-2k}\right\rangle_{\substack{g-1\leq\lambda\leq 2(g-1)\\\lambda-g+1\leq k\leq\lfloor\frac{\lambda}{2}\rfloor}}&
\Bigg{\rangle}.
\end{align*}
Here, the blocks within each megablock are ordered first by decreasing $\lambda$, then by decreasing $k$.
\end{thm}
We proceed by analogy with \cite{tevelevbraid}*{Section 5}, beginning with the reordering trick. This works for any $d\leq 2g$ with $i_d\leq\floor{\frac{d-1}{2}}$.
\begin{lem}[cf. \cite{tevelevbraid}*{Theorem 5.3}] \label{lem:sod_reordered}
    With notation as in \cref{thm:sod}, we have the following semi-orthogonal decomposition of $D^b(M_{i_d}(d))$:
    \begin{equation}\label{eq:sod_reordered1}
    \left\langle 
    \left\langle \bLambda^{-k}\boxpow{\cF^\dual}{\lambda-2k}\right\rangle_{\substack{\lambda-k\leq i_d-m_2\\\lambda-2k,k\geq 0}}, \left\langle T_1\bLambda^{-k}\boxpow{\cF^\dual}{\lambda-2k}\right\rangle_{\substack{\lambda-k\leq i_d-m_1\\\lambda-2k,k\geq 0}},\left\langle T_2\bLambda^{-k}\boxpow{\cF^\dual}{\lambda-2k}\right\rangle_{\substack{\lambda-k\leq i_d\\\lambda-2k,k\geq 0}}
    \right\rangle.
    \end{equation}
    The blocks within each megablock are ordered first by decreasing $\lambda$, then by decreasing $k$.
\end{lem}
\begin{proof}
The blocks in each megablock are the same as those in \eqref{eq:sod1} with $\lambda=j+2k$. As they are already ordered by decreasing $k$, it suffices to show that we can move blocks with smaller $\lambda$ to the right of blocks with larger $\lambda$, i.e., $\langle\bLambda^{-k}\boxpow{\cF^\dual}{\lambda-2k}\rangle\subset {}^\perp\langle\bLambda^{-k'}\boxpow{\cF^\dual}{\lambda'-2k'}\rangle$ for $\lambda<\lambda'$. This follows from \cref{lem:bl_vanishing} below.
\end{proof}
\begin{proof}[Proof of \cref{thm:sod_2g}]
When $d=2g$, we have $i_d=g-1$, $m=2$, $T_1=\theta\bLambda$, and $T_2=\theta^2\bLambda$, so \eqref{eq:sod_reordered1} becomes
$$\left\langle \left\langle\bLambda^{-k}\boxpow{\cF^\dual}{\lambda-2k}\right\rangle_{\substack{\lambda-k\leq g-2\\\lambda-2k,k\geq 0}}, \left\langle(\theta\bLambda)\bLambda^{-k}\boxpow{\cF^\dual}{\lambda-2k}\right\rangle_{\substack{\lambda-k\leq g-1\\\lambda-2k,k\geq 0}},\left\langle(\theta^2\bLambda)\bLambda^{-k}\boxpow{\cF^\dual}{\lambda-2k}\right\rangle_{\substack{\lambda-k\leq g-1\\\lambda-2k,k\geq 0}}\right\rangle.$$
We take the part of the third megablock with $\lambda\leq g-2$ and tensor by $\omega_{M_{g-1}(2g)}=\theta^{-3}\bLambda^{-1}$, moving it to the far left. Tensoring everything by $\bLambda^{\left\lfloor \frac{g-2}{2}\right\rfloor}$ proves the theorem.
\end{proof}

\section{Modified Plain Weave}\label{sec:plainweave}

\subsection{Main result} We proceed with the Plain Weave, cf. \cite{tevelevbraid}*{Section 6}. While the spirit of the argument is the same, there are some technical complications in even degree. Notably, not every pair $(F,s)$ with $F$ a semistable bundle is stable, so $M$ parameterizes only an open substack of all such pairs.

\begin{ntn}\label{ntn:stacks}
We denote by $\mc N$ the  stack of rank $2$ semistable bundles on $C$ with determinant $\Lambda$ (and with $\mbb G_m$ as a generic inertia group) and by $\mbb N$ its rigidification (with trivial generic stabilizers). Concretely, we work with the quotient stacks $\mc N=[Q/\mathrm{GL}(\mbb V)]$ and $\mbb N=[Q/\mathrm{PGL}(\mbb V)]$ where $\mbb V$ is a vector space of dimension $2+2m$ for some large $m$ and $Q$ is an appropriate locally closed subscheme of the Quot scheme parameterizing quotients of $\mbb V\otimes \mc O_C(-mp)$ for some fixed point $p\in C$ (see \cite{kosekitodaderived}*{Section~4} for details).  Here $\mc N$ and $\mbb N$ are smooth algebraic stacks and $Q$ is a smooth quasi-projective variety. The generic inertia group of $\mc N$ is identified with the center $\mbb G_m\subset\mathrm{GL}(\mbb V)$. We have morphisms of stacks 
$$\begin{tikzcd}
    M\arrow[r]& \mc N\arrow[r,"\rho"]&\mbb N\arrow[r]& {\mc SU_C(2)},
\end{tikzcd}$$
where $M=M_{g-1}(\Lambda)$ is the moduli space of stable pairs, $M\to \mc N$ is the forgetful morphism, and ${\mc SU_C(2)}$ is the coarse moduli space (of both $\mc N$ and $\mbb N$) as well as the GIT quotient of $Q$ by $\mathrm{PGL}(\mbb V)$. We do not notationally distinguish between universal bundles $\mc F$ on $M\times C$ or $\mc N\times C$, nor those on other spaces that carry them appearing below;
similarly for $\bLambda=\det \mc F_x$ and $\theta$, the (pullback of the) ample generator of $\Pic {\mc SU_C(2)}$ \cite{drezetnarasimhanpicard}. Unlike in the odd degree case, neither $\mc F$ nor any line bundle twist of $\mc F$ descends to $\mbb N$
or any open substack of it ~\cite{ramananmoduli}*{Theorem 2}. 
However,
twisted tensor vector bundles 
    $\bLambda^{k}\boxpow{\cF^\dual}{2k}$
    on  $\mc N\times \Sym^{2k}C$ 
    have weight $0$ with respect to
    $\mbb G_m$
    and therefore descend to 
    $\mbb N\times \Sym^{2k}C$
    for every $k\ge 0$.

Let $\pi_{\mc N}:\mc N\times C\to\mc N$ be the projection. We write $R\pi_{\mc N*}\cF =[\mc A\overset{u}{\to} \mc B]$ for $\mc A$ and $\mc B$ vector bundles on $\mc N$ of ranks $a$ and $b$, respectively, where $a-b=2$ and $\mbb G_m$ acts with weight $1$ on the fibers of both \cite{kosekitodaderived}*{Lemma~4.4}. Let $\alpha:\mc A\to \mc N$, where we use the same notation for vector bundles and their total spaces. Then $u$ gives a section of the vector bundle $\alpha^*\mc B$ over $\mc A$. Let $\mc Z\subset \mc A$ be the vanishing locus of this section, let $\mc A^\circ\subset \mc A$ be the complement of the zero section, and let $\mc Z^\circ=\mc Z\cap \mc A^\circ$, which is the stack of pairs $\{(F,s):F\in \mc N, s\in H^0(F)\smallsetminus\{0\}\}$ (see \cite{kosekitodaderived}*{Lemma 4.5(i)}). We have a diagram
    \begin{equation}\label{eq:diagram}
        \begin{tikzcd}
        &\mc Z\arrow[r,hook]&\mc A\arrow[r,"\alpha"]&\mc N\arrow[d,"\rho"]\\
        M\arrow[r,hook,"j"]& \mc Z^\circ\arrow[r,hook]\arrow[rr,bend right,"\zeta"]\arrow[u,hook]&\mc A^\circ\arrow[u,hook]&\mbb N.
        \end{tikzcd}
    \end{equation}
\end{ntn}

\begin{prp}\label{prp:zetaproperties} In the notation of \eqref{eq:diagram}:
\begin{prplist}
    \item\label{prp:Zcirckoszul} 
    In $D^b(\mc A^\circ)$, $\mc O_{\mc Z^\circ}$ is isomorphic to the Koszul complex
$$
\left[\bigwedge^b\alpha^*\mc B^\dual|_{\mc A^\circ}\to\ldots\to
\alpha^*\mc B^\dual|_{\mc A^\circ}\to\mc O_{\mc A^\circ}
\right].$$
\item\label{prp:pushO} We have $R\zeta_*\mc O_{\mc Z^\circ}=\mc O_{\mbb N}$.
\item\label{prp:relativedualizing} The relative dualizing sheaf for $\zeta$ is $\theta\bLambda^{-1}$.
\end{prplist}
\end{prp}

\begin{proof}
We abuse notation and denote the morphism 
${\mc A^\circ}\to\mbb N$ by $\zeta$.
It is a (twisted) projective bundle with fiber $\mbb P^{a-1}$.  By \cite{kosekitodaderived}*{Lemma 4.5(i)}, the stack $\mc Z^\circ$ is smooth of dimension $3g-2$.
In particular, its codimension in  $\mc A^\circ$ is equal to $b$. The first claim follows. Since 
$R\zeta_*\mc O_{\mc A^\circ}=\mc O_{\mbb N}$
and $R\Gamma(\mbb P^{a-1},\mc O(-k))=0$
for $k=1,\ldots,a-1$, the second claim follows from the first. 
Indeed, since the claim is local on $\mbb N$, we can trivialize~$\mc B$, and subsequently replace
$\alpha^*\mc B^\dual|_{\mc A^\circ}$
with a direct sum of $b$ copies of $\mc O_{\mc A^\circ}(-1)$.
Finally, since 
$\mc A^\circ$ has $\mbb G_m$-weight $1$, the dualizing sheaf for 
$\mc A^\circ\to\mbb N$ is isomorphic to $\alpha^*\det\mc A^\dual$. Therefore, by adjunction and ignoring pullback by $\alpha$,
$$\omega_\zeta\cong\det\mc A^\dual\otimes\det\mc B\cong
(\det\pi_{\mc N!}\mc F)^\dual,$$
which is isomorphic to $\theta\bLambda^{-1}$
by \cite{narasimhanderived}*{Proposition 2.1}.
\end{proof}

In the diagram \eqref{eq:diagram}, $j$ is the inclusion of the GIT-semistable locus with respect to a certain line bundle $L_\ell\otimes \chi_0^\epsilon$ on $\mc Z^\circ$ \cite{kosekitodaderived}*{Section~4}. This gives rise to windows embeddings $D^b(M)\cong\mb G_w\subset D^b(\mc Z^\circ)$ \cite{halpernleistnerwindows}. In the following proposition, we calculate the relevant weights (cf. \cite{tevelevtorresbgmn}*{Lemma 3.17 and Theorem 3.21}).

\begin{prp}\label{prp:Zcircwindows}
    With respect to the semistable locus $M\hookrightarrow \mc Z^\circ$:
    \begin{prplist}
        \item\label{prp:kempfness} There is a unique Kempf--Ness stratum with associated window width $\eta=g$.
        \item\label{prp:weights} Objects in the subcategory $\langle \theta^x\bLambda^y \boxpow{\cF^\dual}{z}\rangle\subset D^b(\mc Z^\circ)$ have weights in the range $[-y,z-y]$.
    \end{prplist}
\end{prp}
\begin{proof}
As in \cite{tevelevtorresbgmn}*{Section~3}, we identify $\mc A^\circ$ with the quotient stack $[X/\mathrm{PGL}(\mbb V)]$, where  $X\subset Q\times\mbb P(\mbb V)$ is 
    a closed subvariety
    parameterizing quotients $\phi:\,\mbb V\otimes \mc O_C(-mp)\to F$ together with a global section 
    $s(mp)\in\mbb P(\mbb V)$. 
    The complement $\mc U=\mc A^\circ\smallsetminus M$ is a closed substack of pairs $(F,s)$ such that $F$ contains a degree $g$ line subbundle $L$ and $s\in H^0(C,L)\smallsetminus\{0\}$. By semistability of $F$, $L=\mc O_C(D)$ is unique, where $D\in\Sym^gC$ is the vanishing locus of $s$. We have a short exact sequence $0\to \mc O_C(D)\to F\to\Lambda(-D)\to 0$. The unstable locus $S\subset X$ is the preimage of~$\mc U$. For points in $S$, the quotient $\phi$ is given by a block upper-triangular matrix corresponding to a splitting $\mbb V=\mbb V_1\oplus \mbb V_2$, where $\mbb V_1\otimes\mc O_C(-mp)=\phi^{-1}(\mc O_C(D))$ and $\dim\mbb V_1=\dim\mbb V_2=r$. Furthermore, we have $s(mp)\in\mbb P(\mbb V_1)$. The destabilizing one-parameter subgroup $\lambda:\,\mbb G_m\to \mathrm{PGL}(\mbb V)$ is given by sending $\lambda(t)= \operatorname{diag}(u,\ldots,u,u{}^{-1},\ldots,u{}^{-1})$, where $u{}^r=t$
    (see the proof of \cite{tevelevtorresbgmn}*{Lemma 3.17}). 
    It~follows that there is only one 
    Kempf--Ness stratum $(\lambda,Z,S)$,
    where $Z=X^\lambda$ is the locus of split quotients $\mbb V_1\otimes\mc O_C(mp)\oplus\mbb V_2\otimes\mc O_C(mp)\to\mc O_C(D)\oplus\Lambda(-D)$ with $s(mp)\in\mbb P(\mbb V_1)$.
    
    Arguing as in the proof of \cite{tevelevtorresbgmn}*{Lemma 3.17}, the window width $\eta=\operatorname{weight}\mc N^*_{S/X}|_Z$ is equal to the codimension of $S$ in~ $X$. Since $\dim\mc A^\circ=3g-2$, to finish the proof of (a), it suffices to show that $\dim\mc U=2g-2$.
    This follows from the description of $\mc U$ above. Indeed, $\Sym^gC$ has dimension~$g$ and the space of extensions $\mathrm{Ext}^1(\Lambda(-D),\mc O_C(D))$
    has dimension $g-1$ when $2D\not\sim\Lambda$ (and $g$ otherwise). Furthermore, points in $\mc U$ have generic stabilizers $\mathbb G_m$ acting trivially on $\mathbb V_1$ and non-trivially on $\mathbb V_2$.
    This gives  $\dim\mc U=g+(g-1)+(-1)=2g-2$ and proves (a). The proof of part (b) is entirely analogous to the proof of 
    \cite{tevelevtorresbgmn}*{Theorem 3.21}
\end{proof}

\begin{cor}\label{cor:windows}
    The blocks in \cref{thm:sod_2g}, with the kernels now regarded as objects on $D^b(\mc Z^\circ\times \Sym^\ell C)$, give a semiorthogonal decomposition of the windows subcategory $\mb G=\mb G_{-\floor{g/2}}\subset D^b(\mc Z^\circ)$.
\end{cor}
\begin{proof}
    Using \cref{prp:weights}, one checks that objects in those blocks have weights in $[-\floor{\frac{g}{2}},g-\floor{\frac{g}{2}})$, so they are contained in $\mb G$. For example, the $(\lambda,k)$ block in the first megablock has weights in the range $$\left[k-\floor{\frac{g-2}{2}},\lambda-k-\floor{\frac{g-2}{2}}\right]\subseteq \left[-\floor{\frac{g-2}{2}},g-2-\floor{\frac{g-2}{2}}\right]=\left[-\floor{\frac{g-2}{2}},g-\floor{\frac{g}{2}}\right).$$
    The windows embedding $D^b(M)\cong\mb G\subset D^b(\mc Z^\circ)$ is right inverse to the restriction 
    $j^*:\,D^b(\mc Z^\circ)\to D^b(M)$, so~$j^*$ gives an equivalence $\mb G\cong D^b(M)$ taking each block in $\mb G\subset D^b(\mc Z^\circ)$ to the corresponding block in~$D^b(M)$. The result follows.
\end{proof}

\begin{ntn}
We introduce full subcategories $\mb K=\{X\in D^b(\mc Z^\circ):\zeta_*X=0\}$ and $\mb K^\dual=\{X^\dual:X\in \mb K\}$.
\end{ntn}
\begin{rem}\label{rem:Kdual}
Note that $\theta\otimes \mb K=\mb K$ by the projection formula and $\mb K^\dual=\bLambda\otimes \mb K$ by coherent duality and \cref{prp:relativedualizing}.
\end{rem}

\begin{thm}[Plain Weave]\label{thm:plainweave}
There is a semiorthogonal decomposition 
\begin{equation}\label{eq:Gsod}
\mb G=\langle\mb L,\mb i,\ii,\iii,\iv,\mb R\rangle
\end{equation}
where $\mb L\subset \mb K$, $\mb R\subset \mb K^\dual$, and
\begin{align*}
\mb i =\left\langle
\theta^{-1}\bLambda^{m}\boxpow{\cF^\dual}{2m}
\right\rangle
_{0\leq m\leq \floor{\frac{g-2}{2}}},&\quad
\ii = \left\langle
\bLambda^{m}\boxpow{\cF^\dual}{2m}
\right\rangle
_{0\leq m\leq \floor{\frac{g-2}{2}}},\\
\iii=\left\langle\theta\bLambda^{m}\boxpow{\cF^\dual}{2m}
\right\rangle
_{0\leq m\leq \floor{\frac{g-1}{2}}},&\quad
\iv=\left\langle\theta^2\bLambda^{m}\boxpow{\cF^\dual}{2m}\right\rangle_{0\leq m\leq \floor{\frac{g-1}{2}}},
\end{align*}
with each megablock is ordered by increasing $m$.
\end{thm}

We postpone the proof of \cref{thm:plainweave} to the next subsection. We can now precisely state and prove our main result.

\begin{thm}\label{thm:main_precise}
    The admissible subcategory 
    $\mc D=\langle\mb i,\ii,\iii,\iv\rangle\subset\mb G$
    is a noncommutative resolution of sin-gularities of $\mc SU_C(2)$. Furthermore, this resolution agrees with the resolution
    of \cite{padurariu}*{Theorem 1.1} (defined there in a more general context of symmetric stacks). This resolution is strongly crepant \cite{kuznetsovcatres} if $g$ is even.
\end{thm}

\begin{proof}
    The blocks in $\mb i$, $\ii$, $\iii$, $\iv$ are exactly those from \cref{thm:sod_2g} which are pulled back from $\mbb N$. Indeed, the bundle $\theta^i\bLambda^j\boxpow{\cF^\dual}{k}$ on $\mc N$ has weight $2j-k$ with respect to the $\mbb G_m$ action, so it descends to $\mbb N$ exactly when $k=2j$.
    By \cref{prp:pushO}, it follows that $\mc D$ is isomorphic (via $R\zeta_*$) to an
    admissible subcategory $\mc D'\subset D^b(\mathbb N)$ given by the same Fourier--Mukai kernels as $\mc D$ (regarded as objects in $D^b(\mbb N\times \Sym^k C)$).
    
    By the Luna slice theorem, analytic-locally near a complex point $p\in\mc SU_C(2)$, the stack $\mathbb N$ is isomorphic to an analytic neighborhood of the origin in the stack $[N/G]$, where $G$ is the stabilizer of a split bundle $F=\mc O_C(D)\oplus\Lambda(-D)$
    from the S-equivalence class of~$p$ and $N$ is its normal bundle in $\mathbb N$. According to \cite{kirwanhomology}, there are two cases. In the first case,  
    $\Lambda\cong\mc O_C(2D)$,  $G=\text{PGL}_2$, and $N=\mathfrak{sl}_2\otimes\mathbb C^g$.
    In the second case, $\Lambda\not\cong\mc O_C(2D)$, 
    $G\cong \mathbb G_m$ is the maximal torus of $\text{PGL}_2$, and $N=\mathbb C^{g-1}\oplus\mathbb C^{g-1}\oplus\mathbb C^g$. A primitive one-parameter subgroup 
    $\lambda:\,\mbb G_m\to G$, given by sending $\lambda(t)= \operatorname{diag}(u,u^{-1})$ where $u^2=t$, has weights $1,-1,0$ on $N$, each with multiplicity $g$ in the first case, and with multiplicities $g-1,g-1,g$ in the second case. It follows that $\mathbb N$ is a symmetric stack 
    satisfying assumptions A, B, C of \cite{padurariu}, 
    so all of its results apply. The window width 
    $\eta=\text{weight}_\lambda\det\mathbb L^{>0}=\text{weight}_\lambda N^{>0}-\text{weight}_\lambda\mathfrak g^{>0}$ is equal to $g-1$ in both cases. By  \cite{padurariu}*{Theorem~1.1}, it follows that the full subcategory $\mc D''=\{X\in D^b(\mathbb N):-\frac{g-1}{2}\le\operatorname{weight}_\lambda X\le \frac{g-1}{2}\}$
    is admissible and provides a noncommutative resolution of singularities of $\mc SU_C(2)$. 
    
    Next, we will show that 
    $\mc D'=\mc D''$. Let $X\in\mc D'$.
        The computation of  
$\text{weight}_\lambda X$ is the same as the computation of   $\text{weight}_\lambda \zeta^*X$ in 
Proposition~\ref{prp:Zcircwindows}. 
By Lemma~\ref{lem:littleweights} below,
the weights of objects in subcategories $\mb i$, $\ii$, $\iii$, $\iv$
are in the interval $[-\frac{g-1}{2}, \frac{g-1}{2}]$.
It follows that $\mc D'\subset \mc D''$.
On the other hand, let $X\in\mc D''$. 
Since $\text{weight}_\lambda X=\text{weight}_\lambda \zeta^*X$ as above,
we have $\zeta^*X\in\mb G$. 
With respect to the semiorthogonal decomposition~\eqref{eq:Gsod} of Theorem~\ref{thm:plainweave}, we have
$\Hom(\zeta^*X,\mb L)=0$ by projection formula and
$\Hom(\mb R, \zeta^*X)=0$ by coherent duality.
It follows that $\zeta^*X\in\mc D$, and therefore $X\in \mc D'$.
We conclude that 
    $\mc D'=\mc D''$.
    
It remains to show that $\mc D$ is a strongly crepant noncommutative resolution of $\mc SU_C(2)$ if $g$ is even. By the discussion above, we can view $\mc D$ as an admissible subcategory 
of $D^b(M)$, $D^b(\mc Z^{0})$, or $D^b(\mbb N)$, where
in every case the pullback functor 
$\mathrm{Perf}(\mc SU_C(2))\to\mc D$ is the usual pullback. Furthermore, it endows $\mc D$ with the structure of an $\mc SU_C(2)$-linear category via the usual tensor product.
We will view $\mc D$ as a subcategory of~$D^b(M)$. Let $f:\,M\to SU_C(2)$ be the forgetful morphism.
For every $B\in\mc D$,
the functor $\mathrm{Perf}(\mc SU_C(2))\to\mc D$,
$A\mapsto f^*A\otimes B$
 has a right adjoint functor $\mc D\to D^b(\mc SU_C(2))$ given by 
$C\mapsto Rf_*\circ R\mc Hom_M(B,C)$.
According to \cite{kuznetsovcatres}, 
in order to show that $\mc D$ is strongly crepant, we need to show that the identity functor on $\mc D$ is a relative Serre functor for $\mc D$ over $\mc SU_C(2)$. That is, we must give a functorial isomorphism
$$
R\mc Hom_{\mc SU_C(2)}(Rf_*\circ R\mc Hom_M(B,C),\mc O_{\mc SU_C(2)})\cong
Rf_*\circ R\mc Hom_M(C,B)
$$
for $B,C\in \mc D$. By coherent duality for $f$, it is in fact sufficient to establish a functorial isomorphism
\begin{equation*}\label{eq:fhomisom}
Rf_*\circ R\mc Hom_M(C,B\otimes\omega_f[1])\cong
Rf_*\circ R\mc Hom_M(C,B),
\end{equation*}
where $\omega_f$ is the dualizing line bundle for $f$.

By Lemma~\ref{lem:mainmutation2},
there exists 
a  morphism $\gamma:\,\mc O\to \theta\otimes\bLambda^{-1}[1]$ 
in $D^b(\mc Z^\circ)$
whose image under $\zeta_*$ is an isomorphism $\mc O\to \mc O$.
We pull back $\gamma$ via the open immersion
$j:\,M\to \mc Z^\circ$
to 
a morphism 
$j^*\gamma:\,\mc O\to \omega_f[1]$
on $M$.
Indeed, $\omega_f\cong
\omega_M\otimes\omega_{\mc SU_C(2)}^{-1}\cong j^*(\theta\otimes\bLambda^{-1})$ by \cite{thaddeusstablepairs}*{Section 5}.
We will show that $j^*\gamma$ induces a stronger  isomorphism
\begin{equation}\label{eq:zetahom}
R(\zeta\circ j)_*\circ R\mc Hom_M(C,B)
\cong
R(\zeta\circ j)_*\circ R\mc Hom_M(C,B\otimes\omega_{f}[1]).
\end{equation}

By Corollary~\ref{cor:windows} 
and Lemma~\ref{lem:littleweights}, the weights of  objects
$B$, $C$, and $B\otimes\omega_{f}[1]$
belong to the interval $[-\lfloor{\frac{g-1}{2}}\rfloor,\lfloor{\frac{g-1}{2}}\rfloor+1]$,
which, if $g$ is even, has width smaller than the window width $\eta=g$ for the immersion $j$.
The GIT construction of $M$ is local over $\mbb N$,
so the quantization theorem 
\cite{halpernleistnerwindows}*{Theorem 3.29}
gives vertical isomorphisms
in the following commutative diagram:
\[\begin{tikzcd}
	{R(\zeta\circ j)_*\circ R\mc Hom_M(C,B)} & {R(\zeta\circ j)_*\circ R\mc Hom_M(C,B\otimes\omega_{f}[1])} \\
	{R\zeta_*\circ R\mc Hom_{\mc Z^\circ}(C,B)} & {R \zeta_*\circ R\mc Hom_{\mc Z^\circ}(C,B\otimes\theta\otimes\bLambda^{-1}[1])}
	\arrow[from=1-1, to=1-2]
	\arrow["\cong"', from=1-1, to=2-1]
	\arrow["\cong", from=1-2, to=2-2]
	\arrow[from=2-1, to=2-2]
\end{tikzcd}\]
The bottom horizontal morphism of this diagram
is an isomorphism
by the projection formula for the morphism
${\mc Z^\circ}\to\mbb N$ and the fact that the morphism 
$\gamma:\,\mc O\to \theta\otimes\bLambda^{-1}[1]$
pushes forward to an isomorphism~$\zeta_*\gamma$.
It follows that the top horizontal morphism 
is also an isomorphism, proving
\eqref{eq:zetahom}.
\end{proof}

\begin{rem}\label{rem:motive}
 For symmetric stacks $\mc X$ satisfying various assumptions, a noncommutative motive $\mb D^{nc}(\mc X)$ is constructed in \cite{padurariu}*{Section 5}. It categorifies the intersection cohomology of the coarse moduli space of~$\mc X$. When $g$ is even, $\mb D^{nc}(\mbb N)$ is the motive of $\mc D$ by its construction in \cite{padurariu}, since then $\frac{g-1}{2}$ is a half-integer.
\end{rem}

\begin{lem}\label{lem:littleweights}
In the setup of Corollary~\ref{cor:windows}, objects in subcategories $\mb i$, $\ii$, $\iii$, $\iv$
have weights in  $[-\frac{g-1}{2}, \frac{g-1}{2}]$. 
\end{lem}

\begin{proof}
    The same calculation as the proof of Corollary~\ref{cor:windows} using Proposition~\ref{prp:Zcircwindows}.
\end{proof}

\subsection{Proof of the Plain Weave} We require several lemmas, to be proved at the end of this section. All mutations take place in $D^b(M)\cong \mb G\subset D^b(\mc Z^\circ)$.


\begin{lem}\label{lem:makebars}
    For $0\leq \lambda\leq 2(g-1)$, there is a mutation
    \begin{equation}\label{eq:makebars}
        \left\langle\bLambda^{-k}\boxpow{\cF^\dual}{\lambda-2k}\right\rangle_{0\leq k\leq\lfloor\frac{\lambda}{2}\rfloor,\lambda-k\leq g-1}\longrightarrow \left\langle\bLambda^{-k}\boxpow{\bar{\cF^\dual}}{\lambda-2k}\right\rangle_{0\leq k\leq\lfloor\frac{\lambda}{2}\rfloor,\lambda-k\leq g-1}
    \end{equation}
    where the blocks are ordered by decreasing $k$ on the left, and by increasing $k$ on the right.
\end{lem}

\begin{lem}\label{lem:pushtozero}
    For $0\leq 2k+1\leq g-1$, we have $\langle\bLambda^k\boxpow{\cF^\dual}{2k+1}\rangle\subset \mb K$ and $\langle\bLambda^{k+1}\boxpow{\cF^\dual}{2k+1}\rangle\subset \mb K^\dual$.
\end{lem}

\begin{lem}[cf. \cite{tevelevbraid}*{Theorem 6.3}]\label{lem:mainmutation1}
For $\ell=0,1$ and $\ell\leq k\leq\frac{g-1}{2}$, let $\mc D^\ell=\mc D^{\ell}_{g-1}\subset D^b(M\times \Sym^\ell C)$ as in \cref{ntn:ft}. Let $\Phi:D^b(\mc Z^\circ)\to {}^\perp\langle \bLambda^k\boxpow{\cF^\dual}{2k}\rangle$ be the semiorthogonal projector. Then $\Phi(\langle\bLambda^{\ell-k}\cD^\ell\rangle)\subset \mb K^\dual$. 
\end{lem}
    

\begin{lem}\label{lem:mainmutation2}
In $D^b(\mc Z^\circ)$, there exists 
a  morphism $\mc O\to \theta\bLambda^{-1}[1]$ whose image under $\zeta_*$ is an isomorphism $\mc O\to \mc O$.
If $\mb T\subset \zeta^*D^b(\mbb N)$ is a full triangulated subcategory such that $\theta\bLambda^{-1}\otimes \mb T\subset {}^\perp \mb T$, there is a mutation
    \begin{equation*}
    \begin{tikzcd}
    \mb T & \theta\bLambda^{-1}\otimes \mb T\\
    \mb X & \mb T
    \arrow[no head, from=1-2, to=2-1]
    \arrow[no head, crossing over, from=1-1, to=2-2]
    \end{tikzcd}
    \end{equation*}
    where $\mb X\subset \mb K$.
\end{lem}

\begin{figure}[hb]
    \centering
    \includegraphics[width=.9\linewidth]{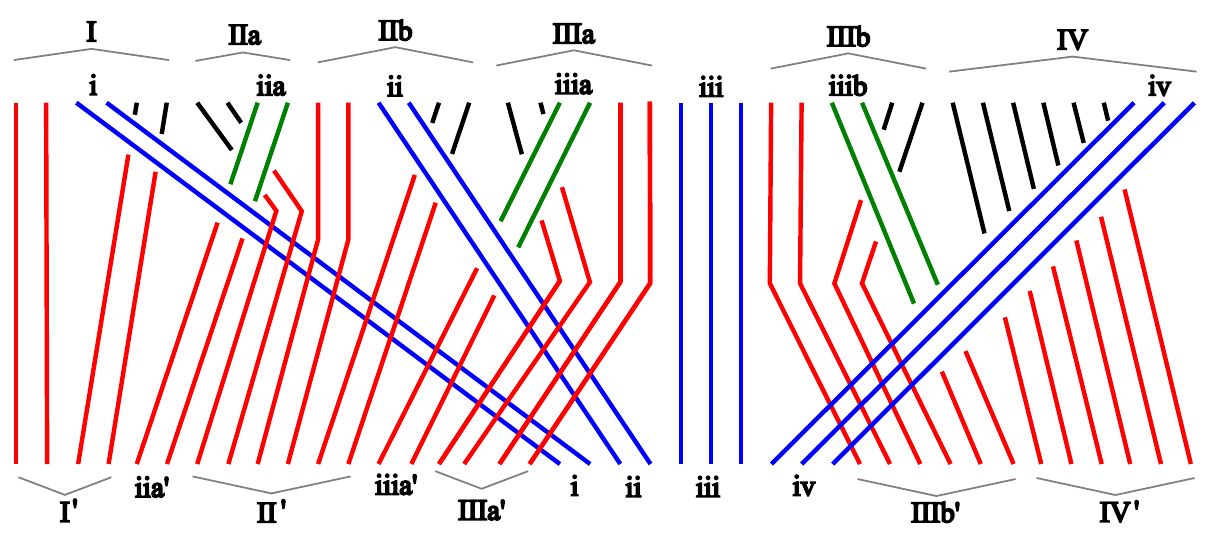}
    \caption{Modified Plain Weave in genus $5$, cf. \cite{tevelevbraid}*{Figure 13}.
    }
    \label{fig:genus5plainweave}
\end{figure}

\begin{proof}[Proof of \cref{thm:plainweave}]
    Denote by $\I,\II,\III,\IV$ the megablocks in \cref{thm:sod_2g}, regarded as subcategories of $\mb G\subset D^b(\mc Z^\circ)$ as in \cref{cor:windows}. We proceed as illustrated in \cref{fig:genus5plainweave}: we take the blocks in $\mb i,\ii,\iii,\iv,$ and move them towards the center, making sure that all blocks in the way mutate into either $\mb K$ on the left or $\mb K^\dual$ on the right. 

    We begin with 
    $$\IV=\left\langle\theta^2\bLambda^{\left\lfloor \frac{g}{2}\right\rfloor-k}\boxpow{\cF^\dual}{\lambda-2k}\right\rangle_{\substack{g-1\leq\lambda\leq 2(g-1)\\\lambda-g+1\leq k\leq\floor{\frac{\lambda}{2}}}}.$$
    The blocks with $\lambda=2\floor{\frac{g}{2}}$ form megablock $\iv$ from the statement (with $m=\floor{\frac{g}{2}}-k$). If $g$ is odd, $2\floor{\frac{g}{2}}=g-1$, so these are the rightmost blocks in $\IV$ as in Figure~\ref{fig:genus5plainweave}; if $g$ is even, the blocks with $\lambda=g-1$ to the right of $\iv$ already lie in $\mb K^\dual$ by \cref{lem:pushtozero}. 
    \begin{clm}\label{clm:mutation_iv}
        For $2\floor{\frac{g}{2}}<\lambda\leq 2(g-1)$ and $\lambda-g+1\leq k\leq \floor{\frac{\lambda}{2}}$, there is a mutation 
        $$\begin{tikzcd}
        \theta^2\bLambda^{\left\lfloor \frac{g}{2}\right\rfloor-k}\boxpow{\cF^\dual}{\lambda-2k} & \iv \\
        \iv & \mb X
        \arrow[no head, from=1-1, to=2-2]
        \arrow[no head, crossing over, from=1-2, to=2-1]
        \end{tikzcd}$$ 
        where $\mb X\subset \mb K^\dual$.
    \end{clm}
    \begin{proof}
        Write $\mb A=\langle\theta^2\bLambda^{\floor{\frac{g}{2}}-k}\boxpow{\cF^\dual}{\lambda-2k}\rangle$. Suppose first that $\lambda=2k$, so $\mb A=\langle\theta^2\bLambda^{\left\lfloor \frac{g}{2}\right\rfloor-k}\rangle$. Take the block $\mb B=\langle\theta^2\bLambda^{m}\boxpow{\cF^\dual}{2m}\rangle$ from $\iv$ with $m=k-\floor{\frac{g}{2}}$ and write $\langle\iv\rangle=\langle \mb B,\mb B'\rangle$. (Note that $0\leq m\leq \floor{\frac{g-1}{2}}$ since $\floor{\frac{g}{2}}< k\leq g-1$.) By \cref{lem:mainmutation1}, we mutate $\langle \mb A,\mb B,\mb B'\rangle \to \langle \mb B,\mb A',\mb B'\rangle$ where $\mb A'\subset \mb K^\dual$. Since $\mb B'\subset \zeta^*D^b(\mbb N)$, we see that $\mb A'$ and $\mb B'$ are fully orthogonal, so we move $\langle \mb B,\mb A',\mb B'\rangle\to\langle \mb B,\mb B',\mb A'\rangle=\langle\iv,\mb A'\rangle$ without any further mutation.

        Now suppose $2k<\lambda$. Let $\mb B=\langle\theta^2\bLambda^{m-1}\boxpow{\cF^\dual}{2m-2},\theta^2\bLambda^{m}\boxpow{\cF^\dual}{2m}\rangle \subset \langle\iv\rangle$, where $m=\lambda-k-\floor{\frac{g}{2}}$. (We have $m>\frac{\lambda}{2}-\floor{\frac{g}{2}}\geq 0$, so $m-1\geq 0$, and $m\leq g-1-\floor{\frac{g}{2}}\leq \floor{\frac{g-1}{2}}$.) Denote by $\boxpow{\cF^\bullet}{\ell}_{tr}=[\boxpow{\cF^\dual}{\ell}\to \mathrm{Ker}^{1-\ell}\boxpow{\cF^\bullet}{\ell}]$ the two-step smart truncation of the complex $\boxpow{\cF^\bullet}{\ell}$ (see \cref{ntn:Fbullet} below). By \cite{tevelevbraid}*{Lemma 4.9, Corollary 4.10}, for every $X\in \mb A$ we have exact triangles $K\to Y\to X \to $ and $H\to Y \to H'\to$ where $Y\in \langle \theta^2\bLambda^{\floor{\frac{g}{2}}-k}\boxpow{\cF^\bullet}{\lambda-2k}_{tr}\rangle$, $K\in \langle\theta^2\bLambda^{\floor{\frac{g}{2}}-k-(\lambda-2k-1)}\rangle=\langle\theta^2\bLambda^{1-m}\rangle$, $H\in\langle\theta^2\bLambda^{-m}\rangle$, and $H'\in \langle\theta^2\bLambda^{1-m}\cD^1\rangle$. It follows from \cref{lem:mainmutation1} that the images of $K,H$ and $H'$ under the semiorthogonal projector onto ${}^\perp \mb B$ lie in $\mb K^\dual$, so the same is true for~$X$. Writing $\langle\iv\rangle=\langle \mb B,\mb B'\rangle$ as above, this gives mutations $\langle \mb A,\mb B,\mb B'\rangle \to \langle \mb B,\mb A',\mb B'\rangle\to\langle \mb B,\mb B',\mb A'\rangle=\langle\iv,\mb A'\rangle$ with $\mb A'\subset \mb K^\dual$, as claimed.
    \end{proof}
    By applying this mutation to each block in $\IV$ with $g\leq\lambda\leq 2(g-1)$ in sequence from right to left, we obtain $\langle\IV\rangle=\langle\iv,\IV'\rangle$ where $\IV'\subset \mb K^\dual$.

    The megablock $$\I=\left\langle\theta^{-1}\bLambda^{\left\lfloor \frac{g-2}{2}\right\rfloor-k}\boxpow{\cF^\dual}{\lambda-2k}\right\rangle_{\substack{0\leq\lambda\leq g-2\\0\leq k\leq\floor{\frac{\lambda}{2}}}}$$ on the left is treated similarly, but the process is mirrored. Megablock $\mb i$ appears as the blocks with $\lambda=2\floor{\frac{g-2}{2}}$. If $g$ is odd, there are additional blocks with $\lambda=g-2$ to the left of the megablock $\mb i$, but they already lie in $\mb K$ by \cref{lem:pushtozero}, so we ignore them as above. We apply the mutation of \cref{lem:makebars} to all blocks in $\I$; the blocks are still ordered by decreasing $\lambda$, but now by increasing $k$ with $\boxpow{\bar{\cF^\dual}}{\lambda-2k}$ in place of~$\boxpow{\cF^\dual}{\lambda-2k}$. 
    \begin{clm}\label{clm:mutation_i}
        For $0\leq\lambda<2\floor{\frac{g-2}{2}}$ and $0\leq k\leq \floor{\frac{\lambda}{2}}$, there is a mutation 
        $$\begin{tikzcd}
        \mb i & \theta^{-1}\bLambda^{\left\lfloor \frac{g-2}{2}\right\rfloor-k}\boxpow{\bar{\cF^\dual}}{\lambda-2k} \\
        \mb X & \mb i
        \arrow[no head, from=1-2, to=2-1]
        \arrow[no head, crossing over, from=1-1, to=2-2]
        \end{tikzcd}$$
        where $\mb X\subset \mb K$.
    \end{clm}
    \begin{proof}
    We prove instead the dual mutation. By \cite{tevelevbraid}*{Lemma 3.4}, we have 
    \begin{equation}\label{eq:dualbar}
        \langle \theta^x \bLambda^y\boxpow{\bar{\cF^\dual}}{z}\rangle^\dual=\langle \theta^{-x} \bLambda^{-y}(\boxpow{\bar{\cF^\dual}}{z})^\dual\rangle=\langle \theta^{-x} \bLambda^{z-y}\boxpow{\cF^\dual}{z}\rangle
    \end{equation}
    Hence it suffices to give a mutation $\langle \mb A,\mb i^\dual\rangle \to \langle\mb i^\dual,\mb A'\rangle$ where $\mb A=\langle\theta\bLambda^{\lambda-k-\floor{\frac{g-2}{2}}}\boxpow{\cF^\dual}{\lambda-2k}\rangle$, $\mb A'\subset \mb K^\dual$, and $\mb i^\dual=\langle\theta\bLambda^m\boxpow{\cF^\dual}{2m}\rangle_{0\leq m\leq \floor{\frac{g-2}{2}}}$ ordered by decreasing $m$. 
    
    From here, the proof is the same as for \cref{clm:mutation_iv}. If $\lambda=2k$, we let $\mb B=\langle\theta\bLambda^m\boxpow{\cF^\dual}{2m}\rangle$ with $m=\floor{\frac{g-2}{2}}-k$ and proceed as above. We need only check that $0\leq m\leq \floor{\frac{g-2}{2}}$, which is clear. Similarly, if $2k<\lambda$, it suffices to check that $1\leq m\leq \floor{\frac{g-2}{2}}$ where $m=(\lambda-2k)-\left(\lambda-k-\floor{\frac{g-2}{2}}\right)=\floor{\frac{g-2}{2}}-k$, which is again clear. 
    \end{proof}
    We thus obtain $\I=\langle \mb I',\mb i\rangle$ where $\mb I'\subset \mb K$. Next, let 
    $$\II_a=\left\langle\bLambda^{\left\lfloor \frac{g-2}{2}\right\rfloor-k}\boxpow{\cF^\dual}{\lambda-2k}\right\rangle_{\substack{g-1\leq\lambda\leq 2(g-2)\\\lambda-g+2\leq k\leq\lfloor\frac{\lambda}{2}\rfloor}},\quad \II_b=\left\langle\bLambda^{\left\lfloor \frac{g-2}{2}\right\rfloor-k}\boxpow{\cF^\dual}{\lambda-2k}\right\rangle_{\substack{0\leq\lambda\leq g-2\\0\leq k\leq\lfloor\frac{\lambda}{2}\rfloor}},$$
    so $\II=\langle\II_a,\II_b\rangle$. Then $\II_b=\theta\otimes \I$ and $\ii=\theta\otimes \mb i$, so $\II_b=\langle \II_b',\ii\rangle$ where $\II_b'=\theta\otimes \I'\subset \mb K$. On the other hand, let $$\ii_a=\left\langle\bLambda^{m-1}\boxpow{\cF^\dual}{2m}\right\rangle_{0\leq m\leq\floor{\frac{g-3}{2}}}$$
    ordered by increasing $m$, which are the blocks in $\II_a$ with $\lambda=2\floor{\frac{g}{2}}$ (where $m=\floor{\frac{g}{2}}-k$). If $g$ is even, there are blocks in $\II_a$ with $\lambda=g-1$ to the right of $\ii_a$, but they are contained in $\mb K$ by \cref{lem:pushtozero}.
    \begin{clm}\label{clm:mutation_ii}
        For $2\floor{\frac{g}{2}}<\lambda\leq 2(g-2)$ and $\lambda-g+2\leq k\leq \floor{\frac{\lambda}{2}}$, there is a mutation 
        $$\begin{tikzcd}
        \bLambda^{\left\lfloor \frac{g-2}{2}\right\rfloor-k}\boxpow{\cF^\dual}{\lambda-2k} & \ii_a \\
        \ii_a & \mb X
        \arrow[no head, from=1-1, to=2-2]
        \arrow[no head, crossing over, from=1-2, to=2-1]
        \end{tikzcd}$$
        where $\mb X\subset \mb K$.
    \end{clm}
    \begin{proof}
        The proof is the same as \cref{clm:mutation_iv}, but with everything tensored with $\bLambda^{-1}\otimes \theta^{-2}$ and with a smaller range of $\lambda$, $k$, and $m$. We need only check that $0\leq m\leq \floor{\frac{g-3}{2}}$ with $m>0$ if $2k<\lambda$, where $m=\lambda-k-\floor{\frac{g}{2}}$. Indeed, we have $0\leq \frac{\lambda}{2}-\floor{\frac{g}{2}}\leq m\leq g-2-\floor{\frac{g}{2}}\leq \floor{\frac{g-3}{2}}$ with $\frac{\lambda}{2}-\floor{\frac{g}{2}}< m$ if $2k<\lambda$. 
    \end{proof}
    Hence we can write $\II_a=\langle \ii_a,\II_a'\rangle$ with $\II_a'\subset\mb K$. Combining with $\II_b$ gives $\II=\langle\ii_a,\II',\ii\rangle$ where $\II'=\langle\II_a',\II_b'\rangle\subset\mb K$. 
    
    Next, we have $\III=\langle\III_a,\iii,\III_b\rangle$, where 
    $$\III_a=\left\langle\theta\bLambda^{\left\lfloor \frac{g}{2}\right\rfloor-k}\boxpow{\cF^\dual}{\lambda-2k}\right\rangle_{\substack{2\floor{\frac{g}{2}}<\lambda\leq 2(g-1)\\\lambda-g+1\leq k\leq\lfloor\frac{\lambda}{2}\rfloor}},
    \quad
    \III_b=\left\langle\theta\bLambda^{\left\lfloor \frac{g}{2}\right\rfloor-k}\boxpow{\cF^\dual}{\lambda-2k}\right\rangle_{\substack{0\leq\lambda<2\floor{\frac{g}{2}}\\0\leq k\leq\lfloor\frac{\lambda}{2}\rfloor}}.$$
    Let $\iii_a=\left\langle\theta\bLambda^{m+1}\boxpow{\cF^\dual}{2m}\right\rangle_{0\leq m\leq \floor{\frac{g-3}{2}}}$ be the blocks with $\lambda=2\floor{\frac{g+2}{2}}$. The blocks in $\III_a$ to the right of $\iii_a$ have $\lambda=2\floor{\frac{g}{2}}+1$ and lie in $\mb K$ already by \cref{lem:pushtozero}. The blocks to the left of $\iii_a$ are processed exactly as with $\II_a$ above (make the substitution $\lambda'=\lambda-2$, $k'=k-1$ and use \cref{clm:mutation_ii}). Hence $\III_a=\langle \iii_a,\III_a'\rangle$ with $\III_a'\subset \mb K$. Similarly, let $\iii_b=\left\langle\theta\bLambda^{m-1}\boxpow{\cF^\dual}{2m}\right\rangle_{0\leq m\leq \floor{\frac{g-2}{2}}}$ be the blocks from $\III_b$ with $\lambda=2\floor{\frac{g-2}{2}}$. The blocks in $\III_b$ with $\lambda=2\floor{\frac{g}{2}}-1$ lie in $\mb K^\dual$; the other blocks are exactly the blocks in $\I$ with $\lambda\leq 2\floor{g-2}{2}$, tensored by $\theta^2\bLambda$. Then \cref{clm:mutation_i} allows us to write $\III_b=\langle\III_b',\iii_b\rangle$ with $\III_b'\subset \mb K^\dual$.

    To summarize, we have a semiorthogonal decomposition $$\mb G=\left\langle\I',\mb i,\ii_a,\II',\ii,\iii_a,\III_a',\iii,\III_b',\iii_b,\iv,\IV'\right\rangle$$
    where $\I',\II',\III_a'\subset\mb K$ and $\III_b',\IV'\subset \mb K^\dual$. It remains to mutate $\ii_a$, $\iii_a$, and $\iii_b$. Observe that the blocks in $\ii_a$ are exactly the blocks from $\mb i$ tensored by $\theta\bLambda^{-1}$. One by one, we mutate each block $\theta\bLambda^{-1}\otimes \mb T$ from $\ii_a$ past the corresponding block from $\mb T$ in $\ii$ using \cref{lem:mainmutation2}; the resulting block $\mb T'\subset \mb K$ is both left and right orthogonal to $\mb T^{\perp}\cap \mb i$, so we may move it to the left of $\mb i$ without further mutation. Hence we obtain $\langle\mb i,\ii_a\rangle=\langle\ii_a',\mb i\rangle$ with $\ii_a'\subset \mb K$. Similarly, $\iii_a=\theta\bLambda^{-1}\otimes\ii$, so $\langle\ii,\iii_a\rangle=\langle\iii_a',\ii\rangle$ with $\iii_a'\subset \mb K$. On the other side, each block in $\iii_b$ is a block from $\iv$ tensored with $\theta^{-1}\bLambda$. The dual of \cref{lem:mainmutation2} allows us to mutate $\langle \theta^{-1}\bLambda\otimes \mb T,\mb T\rangle\to \langle\mb T,\mb T'\rangle$ with $\mb T'\subset \mb K^\dual$, which gives $\langle\iii_b,\iv\rangle=\langle\iv, \iii_b'\rangle$ with $\iii_b'\subset\mb K^\dual$. 
    
    Put together, we have $$\mb G=\left\langle\I',\ii_a',\mb i,\II',\iii_a',\ii,\III_a',\iii,\III_b',\iv,\iii_b',\IV'\right\rangle,$$ where all primed subcategories to the left (resp. right) of $\iii$ lie in $\mb K$ (resp. $\mb K^\dual$). It follows that $\ii$ and $\III_a'$ are both left and right orthogonal; similarly, $\mb i$ is orthogonal to $\II'$, $\iii_a'$, and $\III_a'$, while $\iv$ is orthogonal to $\III_b'$. Thus we may move $\mb i$, $\ii$, and $\iv$ to the center without any mutations, giving $\mb G=\left\langle\mb L,\mb i,\ii,\iii,\iv,\mb R\right\rangle$ where $\mb L=\langle\mb I',\ii_a',\II',\iii_a',\III_a'\rangle\subset \mb K$ and $\mb R=\langle\III_b',\iii_b',\IV'\rangle\subset \mb K^\dual$. This completes the proof.
\end{proof}

\begin{proof}[Proof of \cref{lem:makebars}]
\newcommand{\Abar}{\bar{\mb A}}
\newcommand{\Bbar}{\bar{\mb B}}
\newcommand{\Cbar}{\bar{\mb C}}
Since all blocks on either side of \eqref{eq:makebars} lie in $\mb G$, it suffices to perform the mutation in $D^b(M)$. We begin with the semiorthogonal decomposition of $D^b(M)$ from \cref{lem:sod_reordered}:
\begin{equation} \label{eq:makebarssod1}
    \left\langle \left\langle\bLambda^{-k}\boxpow{\cF^\dual}{\lambda-2k}\right\rangle_{\substack{\lambda-k\leq g-2\\\lambda-2k,k\geq 0}}, \left\langle \theta\bLambda^{1-k}\boxpow{\cF^\dual}{\lambda-2k}\right\rangle_{\substack{\lambda-k\leq g-1\\\lambda-2k,k\geq 0}},\left\langle \theta^2\bLambda^{1-k}\boxpow{\cF^\dual}{\lambda-2k}\right\rangle_{\substack{\lambda-k\leq g-1\\\lambda-2k,k\geq 0}}\right\rangle,
\end{equation}
ordered by decreasing $\lambda$, then decreasing $k$. Using \eqref{eq:dualbar}, we obtain a dual semiorthogonal decomposition:
$$\left\langle \left\langle\theta^{-2}\bLambda^{\lambda-k-1}\boxpow{\bar{\cF^\dual}}{\lambda-2k}\right\rangle_{\substack{\lambda-k\leq g-1\\\lambda-2k,k\geq 0}}, \left\langle \theta^{-1}\bLambda^{\lambda-k-1}\boxpow{\bar{\cF^\dual}}{\lambda-2k}\right\rangle_{\substack{\lambda-k\leq g-1\\\lambda-2k,k\geq 0}},\left\langle \bLambda^{\lambda-k}\boxpow{\bar{\cF^\dual}}{\lambda-2k}\right\rangle_{\substack{\lambda-k\leq g-2\\\lambda-2k,k\geq 0}}\right\rangle,$$
ordered by increasing $\lambda$, then increasing $k$. We tensor the rightmost megablock by $\omega_M=\bLambda^{-1}\theta^{-3}$, moving it to the left, then tensor everything by $\theta^3\bLambda^{3-g}$ to obtain:
$$\left\langle
\left\langle \bLambda^{\lambda-k-g+2}\boxpow{\bar{\cF^\dual}}{\lambda-2k}\right\rangle_{\substack{\lambda-k\leq g-2\\\lambda-2k,k\geq 0}}
\left\langle\theta\bLambda^{\lambda-k-g+2}\boxpow{\bar{\cF^\dual}}{\lambda-2k}\right\rangle_{\substack{\lambda-k\leq g-1\\\lambda-2k,k\geq 0}}, \left\langle \theta^2\bLambda^{\lambda-k-g+2}\boxpow{\bar{\cF^\dual}}{\lambda-2k}\right\rangle_{\substack{\lambda-k\leq g-1\\\lambda-2k,k\geq 0}}\right\rangle,$$
ordered by decreasing $\lambda$, then increasing $k$. Finally, we make the change of variables $\lambda'=2(g-2)-\lambda$, $k'=k-\lambda+g-2$ in the first megablock and $\lambda'=2(g-1)-\lambda$, $k'=k-\lambda+g-1$ in the others to obtain 
\begin{equation}\label{eq:makebarssod2}
    \left\langle \left\langle\bLambda^{-k'}\boxpow{\bar{\cF^\dual}}{\lambda'-2k'}\right\rangle_{\substack{\lambda'-k'\leq g-2\\\lambda'-2k',k'\geq 0}}, \left\langle \theta\bLambda^{1-k'}\boxpow{\bar{\cF^\dual}}{\lambda'-2k'}\right\rangle_{\substack{\lambda'-k'\leq g-1\\\lambda'-2k',k'\geq 0}},\left\langle \theta^2\bLambda^{1-k'}\boxpow{\bar{\cF^\dual}}{\lambda'-2k'}\right\rangle_{\substack{\lambda'-k'\leq g-1\\\lambda'-2k',k'\geq 0}}\right\rangle
\end{equation}
ordered by decreasing $\lambda'$, then increasing $k'$. Denote by $\mb A$, $\mb B$, $\mb C$ the megablocks of $\eqref{eq:makebarssod1}$, and $\Abar$, $\bar{\mb B}$, $\bar{\mb C}$ the megablocks of \eqref{eq:makebarssod2}.
\begin{clm}
    We have $\mb B=\Bbar$.
\end{clm}
\begin{proof}
Since \eqref{eq:makebarssod1} and \eqref{eq:makebarssod2} are both full semiorthogonal decompostions of $D^b(M)$, it suffices to show that $\mb B\subset {}^\perp \Abar\cap \Cbar ^\perp$ and $\Bbar\subset {}^\perp \mb A\cap \mb B^\perp$. By \cite{tevelevtorresbgmn}*{Lemma 10.2}, it suffices to check that 
\begin{equation}\label{eq:makebarsvanishing1}
    R\Hom\left(\theta\bLambda^{1-k'}\boxpow{\bar{\cF^\dual}}{\lambda'-2k'}_{D'},\bLambda^{1-k}\boxpow{\cF^\dual}{\lambda-2k}_{D}\right)=R\Hom\left(\theta\bLambda^{1-k}\boxpow{\cF^\dual}{\lambda-2k}_{D},\bLambda^{1-k'}\boxpow{\bar{\cF^\dual}}{\lambda'-2k'}_{D'}\right)=0
\end{equation}
for any $D\in \Sym^{\lambda-2k}C$, $D'\in \Sym^{\lambda'-2k'}C$, where $\lambda,k$ (likewise $\lambda',k'$) satisfy $k\geq 0$,
$\lambda-2k\geq 0$, and $\lambda-k\leq g-1$. Using $\boxpow{\cF^\dual}{\lambda-2k}_D=\bLambda^{2k-\lambda}\boxpow{\cF}{\lambda-2k}$ and $\boxpow{\bar{\cF^\dual}}{\lambda-2k}_D=\bLambda^{2k-\lambda}\boxpow{\bar\cF}{\lambda-2k}$, \eqref{eq:makebarsvanishing1} becomes
$$ R\Gamma\left(\theta^{-1}\bLambda^{\lambda'-k'-\lambda+k}(\boxpow{\bar\cF}{\lambda'-2k'}_{D'})^\dual\boxpow{\cF}{\lambda-2k}_D\right)=R\Gamma\left(\theta^{-1}\bLambda^{\lambda-k-\lambda'+k'}(\boxpow{\cF}{\lambda-2k}_{D})^\dual\boxpow{\bar\cF}{\lambda'-2k'}_{D'}\right)=0.$$
Recalling that $\theta=\mc O(1,g-1)$, this follows from \cite{tevelevtorresbgmn}*{Theorem 4.1, Remark 4.2} once we verify that $\lambda-2k-g<\lambda-k-(\lambda'-k')<g-\lambda+2k'$. Indeed, $-k-g<-(g-1)\leq -(\lambda'-k')$ and $\lambda-k\leq g-1<g+k'$.
\end{proof}
Hence $\theta^{-1}\bLambda^{-1}\mb B=\theta^{-1}\bLambda^{-1}\Bbar$. Let for $0\leq \lambda\leq 2(g-1)$, let $\mb B_\lambda=\left\langle\bLambda^{-k}\boxpow{\cF^\dual}{\lambda-2k}\right\rangle_{\max(0,\lambda-g+1)\leq k\leq \floor{\lambda/2}}$ (ordered by decreasing $k$), so $\theta^{-1}\bLambda^{-1}\mb B=\langle \mb B_{2(g-1)},\dots,\mb B_0\rangle$; likewise, $\theta^{-1}\bLambda^{-1}\Bbar=\langle \Bbar_{2(g-1)},\dots, \Bbar_0\rangle$.
\begin{clm}
    For $0\leq \lambda\leq 2(g-1)$, we have $\mb B_\lambda=\Bbar_{\lambda}$.
\end{clm}
\begin{proof}
    It suffices to show that $\mb B_\lambda\subset \Bbar_{\lambda'}^\perp$ for $\lambda<\lambda'$ and $\mb B_\lambda\subset {
    }^\perp\Bbar_{\lambda'}$ for $\lambda>\lambda'$. For the first, we must check $$R\Gamma\left(\bLambda^{\lambda-k-\lambda'+k'}(\boxpow{\cF}{\lambda-2k}_{D})^\dual\boxpow{\bar\cF}{\lambda'-2k'}_{D'}\right)=0$$
for any $\max(0,\lambda-g+1)\leq k\leq \floor{\lambda/2}$, $\max(0,\lambda'-g+1)\leq k'\leq \floor{\lambda'/2}$, $D\in \Sym^{\lambda-2k}C$, $D'\in \Sym^{\lambda'-2k'}C$. This follows from \cite{tevelevtorresbgmn}*{Lemma 5.4, Remark 5.7}: we have $\lambda-2k,\lambda'-2k'\leq 2g-1$ and $2(\lambda-k-\lambda'+k')<\lambda-2k-\lambda'+2k'$ immediately, and $\lambda-2k-g<\lambda-k-\lambda'+k'<g-\lambda+2k'+1$ as in the preceding claim. The second containment is proved analogously.
\end{proof}
It remains to show that the blocks of $\mb B_\lambda$ and $\Bbar_{\lambda}$ are related by a mutation. If $\lambda>g-1$, then $\mb B_\lambda=\bLambda^{g-1-\lambda}\mb B_{2(g-1)-\lambda}$, so we may assume $0\leq \lambda\leq g-1$. Proceed by induction on $\lambda$, with $\lambda=0,1$ being trivial. Assume we have the mutation $\mb B_{\lambda-2}\to\Bbar_{\lambda-2}$. Then $\mb B_\lambda=\langle\bLambda^{-1}\mb B_{\lambda-2}, \boxpow{\mc F^\dual}{\lambda}\rangle$ and $\Bbar_\lambda=\langle \boxpow{\bar{\mc F^\dual}}{\lambda},\bLambda^{-1}\Bbar_{\lambda-2}\rangle$. Projecting $\langle\boxpow{\mc F^\dual}{\lambda}\rangle$ onto $(\bLambda^{-1}\mb B_{\lambda-2})^\perp=(\bLambda^{-1}\Bbar_{\lambda-2})^\perp$ and mutating $\mb B_{\lambda-2}\to\Bbar_{\lambda-2}$ completes the proof.
\end{proof}

\begin{proof}[Proof of \cref{lem:pushtozero}]
Let $X\in \langle\bLambda^k\boxpow{\cF^\dual}{2k+1}\rangle$. Then $X$ is a pullback of an object in $D^b(\mc N)$ of weight $-1$ with respect to $\mbb G_m$. 
By \cref{prp:Zcirckoszul}, it suffices to show that $R\zeta_*Y=0$ for every  object 
$Y\in D^b(\mc A^\circ)$
of the form $Y=\bigwedge^k\alpha^*{\mc B}^*\otimes^L\alpha^*{Z}$, where
$k=0,\ldots,b$ and $Z\in D^b(\mc N)$ is an object of weight $-1$. 
Since the claim is local on $\mbb N$, we can replace $Y$  
with $\mc O_{\mc A^\circ}(s)$, where~$s=-1,\ldots,-(b+1)$.
Since $R\Gamma(\mb P^{a-1},\mc O(-s))=0$
for $s=1,\ldots,a-1$ and $b+1=a-1$, the first statement follows. The second follows from \cref{rem:Kdual}.
\end{proof}

\begin{proof}[Proof of \cref{lem:mainmutation1}]
    We mimic the proof of \cite{tevelevbraid}*{Theorem 6.3}. As we will need to work in both $D^b(M)$ and $D^b(\mc Z^\circ)$, we denote the windows embedding by $\iota:D^b(M)\to \mb G\subset D^b(\mc Z^\circ)$. As in \cite{tevelevbraid}*{Lemma 6.7}, it suffices to show that the morphism $$\zeta_*(\bLambda^{-1}\otimes\iota(\bLambda^{\ell-k}\mc O_{M(-D)}))\to \zeta_*(\bLambda^{-1}\otimes\iota\circ\mc P\circ \mc P^L(\bLambda^{\ell-k}\mc O_{M(-D)}))$$ is an isomorphism for any $D\in \Sym^\ell C$, where $\mc P:D^b(\Sym^k C)\to D^b(M)$ is the Fourier--Mukai functor with kernel $\bLambda^{k}\boxpow{\cF^\dual}{2k}$, $\mc P^L$ is its left adjoint (the Fourier--Mukai functor with kernel $(\bLambda^{k}\boxpow{\cF^\dual}{2k})^\dual\otimes \omega_M^\bullet$ \cite{huybrechtsfouriermukai}*{Proposition 5.9}), $M(-D)\subset M$ denotes the locus of stable pairs $(F,s)$ with $s|_D=0$, and the morphism is induced by the unit of adjunction $\mathrm{Id}\implies \mc P\circ \mc P^L$. We compute both sides of this morphism.
    \begin{clm}
        We have $\mc P\circ \mc P^L(\bLambda^{\ell-k}\mc O_{M(-D)})\cong R\pi_{M*}(\bLambda^{k}\boxpow{\cF^\dual(-D)}{2k})[2\ell]$, where $\pi_M:M\times \Sym^k C\to M$ is the projection.
    \end{clm}
    \begin{proof}
        Notice first that by \cite{huybrechtsfouriermukai}*{Corollary 3.40} and \cite{thaddeusstablepairs}*{5.7, 6.1}, we have $$\mc O_{M(-D)}^\dual=\omega_{M(-D)}\otimes \omega_M^{-1}[-2\ell]=\bLambda^\ell\mc O_{M(-D)}[-2\ell]$$ in $D^b(M)$. We have 
        \begin{align*}
            \mc P^L(\bLambda^{\ell-k}\mc O_{M(-D)})
            &=R\pi_{\Sym^k C*}((\bLambda^{2k}\boxpow{\cF^\dual}{2k}\mc O_{M(-D)})^\dual \omega_M^\bullet)[2\ell]\\
            &\cong \left(R\pi_{\Sym^k C*}(\boxpow{\cF(-D)}{2k}|_{M(-D)\times \Sym^k C})\right)^\dual\otimes
       \boxpow{\mc O(-D)}{2k}[2\ell]
        \end{align*} by coherent duality and the projection formula. Since $\cF(-D)|_{M(-D)\times C}$ is the universal family on $M(-D)\times C$, we have $\mc P^L(\bLambda^{\ell-k}\mc O_{M(-D)})\cong\boxpow{\mc O(-D)}{2k}[2\ell]$ by \cite{tevelevtorresbgmn}*{Corollary 7.5}. Applying $\mc P$ proves the claim.
    \end{proof}
    Hence (after shifting by $-2\ell$ for convenience) we have a morphism 
    \begin{equation}\label{eq:morphism}
        \bLambda^{-k}\mc O_{M(-D)}^\dual\to R\pi_{M*}(\bLambda^{k}\boxpow{\cF^\dual(-D)}{2k}),
    \end{equation} which is unique up to scalar as in \cite{tevelevbraid}*{Remark 6.10}. 
    
    It remains to show that applying the functor $\zeta_*(\bLambda^{-1}\otimes\iota(-))$ to \eqref{eq:morphism} yields an isomorphism. Since $R\pi_{\mc Z^\circ*}(\bLambda^{k}\boxpow{\cF^\dual(-D)}{2k})\in D^b(\mc Z^\circ)$ has weights in the range $[-k,k]\subseteq [-\floor{\frac{g}{2}},g-\floor{\frac{g}{2}})$ and restricts via $j^*$ to $R\pi_{M*}(\bLambda^{k}\boxpow{\cF^\dual(-D)}{2k})$, we have $$\iota(R\pi_{M*}(\bLambda^{k}\boxpow{\cF^\dual(-D)}{2k}))=R\pi_{\mc Z^\circ*}(\bLambda^{k}\boxpow{\cF^\dual(-D)}{2k}).$$
    On the other hand, we claim that $\iota(\bLambda^{-k}\mc O_{M(-D)}^\dual)=\bLambda^{-k}\mc O_{\mc Z^\circ(-D)}^\dual$ for $\ell=0,1$, where $\mc Z^\circ(-D)\subset \mc Z^\circ$ denotes the closed substack of pairs $(F,s)$ with $s|_D=0$. If $\ell=0$, this is clear, since $\bLambda^{-k}\in\mb G$. If $\ell=1$, so $D=x\in C$, then $\mc Z^\circ(-x)$ is the codimension-2 vanishing locus of the canonical section of $\mc F_x$. We have a Koszul resolution $\mc O_{\mc Z^\circ(-D)}\cong [\bLambda^{-1}\to \mc F_x^\dual\to \mc O]$, so $\mc O_{\mc Z^\circ(-D)}$ has weights in the range $[0,1]$. Since $[k-1,k]\subseteq [-\floor{\frac{g}{2}},g-\floor{\frac{g}{2}})$ and $j^*\mc O_{\mc Z^\circ(-p)}=\mc O_{M(-p)}$, the claim holds. 
    Thus applying $\iota$ to \eqref{eq:morphism} gives
    \begin{equation}\label{eq:morphism2}
        \bLambda^{-k}\mc O_{\mc Z^\circ(-D)}^\dual\to R\pi_{\mc Z^\circ*}(\bLambda^{k}\boxpow{\cF^\dual(-D)}{2k}),
    \end{equation}
    again unique up to scalar. Moreover, this morphism is not zero: if it were, its cone $\Phi(\bLambda^{-k}\mc O_{\mc Z^\circ(-D)}^\dual)[1]$ would have $R\pi_{\mc Z^\circ*}(\bLambda^{k}\boxpow{\cF^\dual(-D)}{2k})\in \langle \bLambda^{k}\boxpow{\cF}{2k}\rangle$ as a direct summand, which is absurd.
    
    By \cref{prp:relativedualizing}, applying the functor $(\theta \otimes R\zeta_*(\bLambda^{-1}\otimes -))^\dual\cong R\zeta_*((-)^\dual)$ to \eqref{eq:morphism2} gives a morphism
    \begin{equation}\label{eq:morphism3}
        R\zeta_*\left([R\pi_{\mc Z^\circ*}(\bLambda^{k}\boxpow{\cF^\dual(-D)}{2k})]^\dual\right)\to R\zeta_*(\bLambda^k\mc O_{\mc Z^\circ(-D)}),
    \end{equation}
    which we must show is an isomorphism. In fact, it suffices to show that the source and target of \eqref{eq:morphism3} are isomorphic. Indeed, write \eqref{eq:morphism2} as $X\to Y$, where $Y\cong L\zeta^*Z$ (so $R\zeta_*(Y^\dual)\cong Z^\dual$ by the projection formula). Then $\CC\cong \Hom(X,Y)\cong \Hom(L\zeta^*Z^\dual,X^\dual)\cong \Hom(R\zeta_*(Y^\dual),R\zeta_*(X^\dual))$, so \eqref{eq:morphism3} is nonzero and unique up to scalar.

    Recall that $\mc N$ is isomorphic to the moduli stack of rank $2$ vector bundles on $C$ with determinant $\Lambda(-2D)$, with universal family $\cF(-D)$ on $\mc N\times C$. As in \cref{ntn:stacks}, we write $R\pi_{\mc N*}(\cF(-D))=[\mc A'\to \mc B']$, where $\mc A'$ and $\mc B'$ are vector bundles on $\mc N$ with $\mbb G_m$-weight $1$ and ranks $a'$, $b'$, where $a'-b'=2-2\ell$. A polystable vector bundle of the form $\cO(D)\oplus\cO(D)$, where $D$ is an effective divisor of degree $g-l$, has at least a $2$-dimensional space of global sections, so $a'\ge 2$.  Writing $\alpha:\mc Z^\circ\to \mc N$ so that $\zeta=\rho\circ \alpha$, we have 
    \begin{align*}
        R\pi_{\mc Z^\circ*}(\bLambda^{k}\boxpow{\cF^\dual(-D)}{2k})&\cong \bLambda^{-k}\otimes R\pi_{\mc Z^\circ*}(\alpha\times \mathrm{id})^*(\boxpow{\cF(-D)}{2k})\\
        &\cong\alpha^*(\bLambda^{-k}\otimes R\pi_{\mc Z^\circ*}(\boxpow{\cF(-D)}{2k}))\\
        &\cong\alpha^*(\bLambda^{-k}\Sym^{2k}[\mc A'\to \mc B']).
    \end{align*}
    (For the last equality, see the proof of \cite{tevelevbraid}*{Lemma 6.13}). Hence the left hand side of \eqref{eq:morphism3} is the descent 
    of 
    $\bLambda^{k}\Sym^{2k}[\mc A'\to \mc B']^\vee$
    to $\mbb N$ (note that this has $\mbb G_m$-weight $0$). Analogous to \eqref{eq:diagram}, we have a diagram
    $$\begin{tikzcd}
        \mc Z(-D)\arrow[r,hook]&\mc A'\arrow[r,"\alpha'"]&\mc N\arrow[d,"\rho"]\\
         \mc Z^\circ(-D)\arrow[r,hook]\arrow[rr,bend right,"\zeta'"]\arrow[u,hook]&\mc A'^\circ\arrow[u,hook]&\mbb N.
        \end{tikzcd}$$
    The right hand side of \eqref{eq:morphism3} is $R\zeta'_*(\bLambda^k)$. Hence it suffices to prove:
    \begin{clm}\label{clm:rzetalambdak}
        We have $R\zeta'_*(\bLambda^k)\cong \bLambda^{k}\Sym^{2k}[\mc A'\to \mc B']^\vee$.
    \end{clm}
    As in \cref{prp:Zcirckoszul}, we have a Koszul resolution in $D^b(\mc A'^\circ)$: 
    $$\mc O_{\mc Z^\circ(-D)}\cong\left[\bigwedge^{b'}\alpha'^*\mc B'^\dual|_{\mc A'^\circ}\to\ldots\to\alpha'^*\mc B'^\dual|_{\mc A'^\circ}\to\mc O_{\mc A'^\circ}\right].$$
    As above, $\zeta':\mc A'^\circ\to \mbb N$ is a twisted projective bundle with fiber $\PP^{a'-1}$. 
    It follows that 
$$R\zeta_*'(\bLambda^k)\cong
R\zeta_*\left[\bigwedge^{b'}\alpha'^*\mc B'^\dual|_{\mc A'^\circ}\otimes\bLambda^k\to\ldots\to\alpha'^*\mc B'^\dual|_{\mc A'^\circ}\otimes\bLambda^k\to\mc O_{\mc A'^\circ}\otimes\bLambda^k\right].$$
We claim that $R^s\zeta_*'[\bigwedge^{m}\alpha'^*\mc B'^\dual|_{\mc A'^\circ}\otimes\bLambda^k]=0$
for all $s$ when $2k<m\le b'$ and for $s>0$ when $2k\ge m$.
Indeed, we can 
work locally on $\mbb N$, so $\alpha'^*\mc B'$ can be replaced by $\mc O_{\zeta'}(1)^{\oplus b'}$. 
Recall that $\bLambda$ has $\mbb G_m$-weight $2$.
Note that $R\Gamma(\mb P^{a'-1},\mc O(2k-m))=0$ for $2k<m\leq b'$, since then $-a'=2\ell-b'-2< 2k-m<0$. 
It follows that 
$$R\zeta_*'(\bLambda^k)\cong
\left[\zeta_*\bigwedge^{2k}\alpha'^*\mc B'^\dual|_{\mc A'^\circ}\otimes\bLambda^k\to\ldots\to\zeta_*\alpha'^*\mc B'^\dual|_{\mc A'^\circ}\otimes\bLambda^k\to\zeta_*\mc O_{\mc A'^\circ}\otimes\bLambda^k\right]$$
(underived pushforwards, since higher cohomologies vanish).
Since $\alpha'_*\mc O_{\mc A'^\circ}\cong \Sym^\bullet \mc A'^\dual$ (recall that $a'\geq 2)$, computing the zero-weight part gives
$$
R\zeta_*'(\bLambda^k)\cong\left[\bLambda^k\otimes\bigwedge^{2k}\mc B'^\dual\to\ldots\to\bLambda^k\otimes\mc B'^\dual\otimes\Sym^{k-1}\mc A'^\dual\to\bLambda^k\otimes\mc \Sym^{k}\mc A'^\dual\right],$$
which is indeed isomorphic to $
\bLambda^k\Sym^{2k}[\mc A'\to\mc B']^\dual$.

\end{proof}

\begin{proof}[Proof of \cref{lem:mainmutation2}]
    By \cref{prp:pushO,prp:relativedualizing}, coherent duality gives $$R\zeta_*R\mc Hom(\mc O,\theta\otimes\bLambda^{-1}[1])\cong R\mc Hom(\mc O,\mc O)\cong\mc O.$$ Applying $R^0\Gamma$ gives a nonzero morphism $\mc O\to \theta\otimes\bLambda^{-1}[1]$ whose image under $\zeta_*$ is an isomorphism $\mc O\to \mc O$ by construction. We complete this morphism to an exact triangle $\mc O\to \theta\otimes\bLambda^{-1}[1]\to K\to$, where $K\in\mb K$. Tensoring with objects of $\mb T$ gives the required mutation.

\end{proof}

\newcommand{\eN}{\mathcal{N}}
\newcommand{\eF}{\mathcal F}
\newcommand{\etheta}{\theta}

\newcommand{\oN}{\hat{N}}
\newcommand{\oE}{\mathcal{E}}
\newcommand{\oF}{\hat{\mathcal{F}}}
\newcommand{\otheta}{\hat{\theta}}

\section{Quasi-BPS categories and the Hecke Braid}\label{sec:hecke}

\subsection{Review of the Hecke correspondence}
We now relate the semiorthogonal decomposition given in \cite{tevelevbraid}*{Theorem~1.1} (regarding the moduli space of semistable rank $2$ vector bundles with fixed odd determinant) to that of \cref{thm:main_precise} (even determinant) by way of the Hecke correspondence \cite{narasimhanramanangeometry}. We find it convenient to work with slightly different conventions compared to the previous section, which we now describe.

\begin{ntn}\label{ntn:odd}
In this section, we denote by $\eN$ the stack of semistable rank $2$ vector bundles on $C$ with {\it trivial} determinant, with universal family $\eF$ on $C\times \eN$. (This stack is isomorphic to the one in \cref{ntn:stacks}, and their universal families agree up to a twist by a line bundle pulled back from $C$.) The line bundle $\etheta$ on the stack $\eN$ is the pullback 
of the ample generator of the Picard group of its coarse moduli space $\mc SU_C(2)$. We also work with the odd degree counterparts of these objects: fix a point $q\in C$, and let $\oN$ be the moduli space of stable rank $2$ vector bundles on $C$ with determinant $\mc O_C(q)$, with Poincar\'e bundle $\oE$ on $C\times \oN$ normalized such that $\det\oE_q$ is the ample generator $\otheta$ of $\operatorname{Pic}\oN$, where $\oE_q=\oE|_{q\times \oN}$. 
\end{ntn}

We briefly review the Hecke correspondence: $\oE_q$ is a rank $2$ vector bundle on $\oN$ whose associated projective bundle of quotients $P=\mathbb{P}(\oE_q^\dual)$ is a moduli space of stable parabolic bundles on $(C,q)$, i.e., stable rank $2$ bundles $E$ on $C$ with $\det E=\mc O_C(q)$ equipped with a surjection $E\twoheadrightarrow \mc O_q$ onto the skyscraper sheaf at $q$. Indeed, there is a universal family of parabolic bundles $\pi^*\oE\twoheadrightarrow i_*\mc O_\pi(1)$ on $C\times P$, where $\pi:P\to \oN$ and $i:q\times P\hookrightarrow C\times P$, and the map is adjoint to $\oE_q\twoheadrightarrow \mc O_\pi(1)$. (Here and below, we denote morphisms and their base change to $C$ by the same symbol, so $\pi:C\times P\to C\times \oN$ as well.) The kernel of the map $\pi^*\oE\twoheadrightarrow i_*\mc O_\pi(1)$
is a family of rank $2$ bundles of trivial determinant; one can check that this family is semistable, so we get a morphism $\sigma:P\to \eN$ and a short exact sequence
\begin{equation}\label{eq:heckeuniversal}
    0\to \sigma^*\eF \to \pi^*\oE\to i_*\mc O_\pi(1)\to 0,
\end{equation}
on $C\times P$. By \cite{narasimhanderived}*{Proposition 3.3}, we have
\begin{equation}\label{eq:thetaisO1}
    \sigma^*\etheta=\mc O_\pi(1).
\end{equation} Taking determinants in \eqref{eq:heckeuniversal} and pulling back along $i$ yields 
\begin{equation}\label{eq:thetaislambda}
    \sigma^*\bLambda=\pi^*\otheta,
\end{equation} where $\bLambda=\det \eF_q$ as in the previous sections. 

According to \cite{orlovblowup}, we have a semiorthogonal decomposition
\begin{equation}\label{eq:p_sod}
    D^b(P)=\langle \pi^* D^b(\oN),\pi^* D^b(\oN)\otimes \mc O_\pi(1)\rangle.
\end{equation}
By \cite{tevelevbraid}*{Theorem 1.1}, this decomposition has a refinement into blocks equivalent to derived categories of symmetric powers of $C$ (see \eqref{eq:odd_sod} below). The goal of this section, achieved in \cref{thm:hecke_sod_cd}, is to mutate these blocks into blocks pulled back from $\eN$. This is another instance of a (generalized) two-ray game (cf. \cref{cnj:tworay}), as $P$ is a Fano variety with anticanonical bundle $\pi^*\otheta\otimes \mc O_\pi(2)$. To describe the outcome of this mutation, we need twisted analogs of the noncommutative resolution studied in \cref{sec:plainweave} known as quasi-BPS categories \cite{padurariutodaquasibps}.

\subsection{Quasi-BPS categories}
For $w\in \mbb{Z}$, denote by $D^b(\eN)_w\subset D^b(\eN)$ the full subcategory consisting of objects of weight $w$ with respect to the action of $\mbb{G}_m$ by scalar automorphisms. Following \cite{padurariutodaquasibps}*{Definition~3.1}, the \textit{quasi-BPS category} $\mbb{B}_w\subset D^b(\eN)_w$ is defined to be the full subcategory of objects $X$ such that 
\begin{equation}\label{eq:quasibps}
    \operatorname{weight} \nu^* X\subset \left[-\frac{n_\nu}{2},\frac{n_\nu}{2}\right]+\frac{w}{2}\operatorname{weight} \nu^* \bLambda
\end{equation}
for any $\nu:B\mbb{G}_m\to\eN$, where $n_\nu=\operatorname{weight} \det (\nu^* \mathbb{L}_{\eN})^{>0}$. 
\begin{rem}\label{rem:quasibps}
Note that $\mbb{B}_{w+2}=\bLambda\otimes\mbb{B}_{w}$ and $\etheta\otimes\mbb B_w=\mbb B_w$. Also, $\mbb{B}_0$ is exactly the noncommutative resolution $\cD$ of $\mc SU_C(2)$ considered in \cref{sec:plainweave}. 
\end{rem}

\begin{prp}\label{prp:quasibps}
$\sigma^*:\mbb B_w\to D^b(P)$ is a fully faithful admissible embedding for every $w\in \mbb Z$. Moreover, we have $\sigma^*\mbb B_{w+1}\subset {}^\perp\sigma^*\mbb B_{w}$.
\end{prp}
\begin{proof}
    Let $\phi:\mc P\to \eN$ be the map from the total space of the vector bundle  $\eF_q$ with the zero section removed.
    Since $\mathrm{Ext}^1(\mc O_q,E)=\mathrm{Hom}(E,\mc O_q)^\dual=E_q$ for any vector bundle $E$ on $C$, we see that $\mc P$ is isomorphic to the moduli stack of parabolic bundles $E'\to \mc O_q$ whose kernel is semistable of trivial determinant. Hence we have an open immersion $j:P\hookrightarrow \mc P$ as a GIT-semistable locus (see \cite{moonleesymmetric}*{Section 2}), where $\phi\circ j=\sigma$. We claim that both $\phi^*\mbb B_{w}$ and $\phi^*\mbb B_{w+1}$ are contained in the window subcategory $\mathbf{G}=\mathbf{G}_{-\floor{\frac{g-1+w}{2}}}\subset D^b(\mathcal P)$ for this inclusion (see \cite{halpernleistnerwindows}). The argument is very similar to the proof of \cref{thm:main_precise} and \cref{prp:weights}. It follows from \cite{moonleesymmetric}*{Lemma 2.7} that the codimension of the unique stratum of the unstable locus, and hence the corresponding window width, is at least $g$. 
    Furthermore, $n_\nu=g-1$ and $\operatorname{weight}\nu^* \bLambda=-1$ when $\nu$ corresponds to the destabilizing one-parameter subgroup of the stratum. If $g+w=2p$ is even, $\mbb B_w$ and $\mbb B_{w+1}$ have weights in the ranges $[1-p,g-p-1]$ and $[-p,g-p-1]$ respectively (note that the weights must be integers). Similarly, if $g+w=2p+1$, then $\mbb B_{w}$ and $\mbb B_{w+1}$ have weights in the ranges $[-p,g-p-1]$ and $[1-p,g-p-1]$, respectively. In either case, the union of the two ranges is contained in the interval~$[\floor{\frac{g+w-1}{2}},\floor{\frac{g+w-1}{2}}+g)$. This proves the claim.

   
   According to  \cite{orlovblowup}, we have a semiorthogonal decomposition
   $D^b(\mc P)=\langle D^b(\eN)_w,D^b(\eN)_{w+1}\rangle$. 
   It follows that $\phi^*$ gives fully faithful admissible embeddings of $\mbb B_w$ and $\mbb B_{w+1}$ into $\mathbf{G}$ such that $\phi^*\mbb B_{w+1}\subset {}^\perp\phi^*\mbb B_{w}$. Since we have $\sigma^*=j^*\circ \phi^*$ and $j^*|_{\mathbf{G}}$ is fully faithful, this completes the proof.
\end{proof} 
\subsection{The Hecke Braid}
We now describe the promised mutation in $D^b(P)$, which yields semiorthogonal decompositions of these quasi-BPS categories (see \cref{thm:hecke_sod_cd}). For the sake of continuity, the proofs of several lemmas are deferred until the end of the section. We recall \cite{tevelevbraid}*{Theorem 1.1}: we have a semiorthogonal decomposition
\begin{equation}\label{eq:odd_sod}      
    D^b(\oN)=\left\langle \otheta^{-1} \mb A,\otheta^{-1}\mb B,\mb A,\mb B'\right\rangle
\end{equation}
where
\begin{align*}
    \mb A &= \langle A_k\rangle_{0\leq k\leq \floor{\frac{g-2}{2}}}\\
    \mb B &= \langle B_k\rangle_{1\leq k\leq \floor{\frac{g-1}{2}}}\\
    \mb B' &= \langle B_k\rangle_{0\leq k\leq \floor{\frac{g-1}{2}}}\\
    A_k&=\otheta^{2-g+k} \boxpow{\oE}{g-2-2k}\\
    B_k&=\otheta^{2-g+k} \boxpow{\oE}{g-1-2k},
\end{align*}
ordered within each megablock by decreasing $k$ (see \cref{ntn:stablepairs} for tensor bundles). Since $\omega_{\oN}=\otheta^{-2}$, we also have
\begin{equation}\label{eq:odd_sod2}      
    D^b(\oN)=\left\langle\mb A,\mb B',\otheta \mb A,\otheta \mb B\right\rangle.
\end{equation}
Combining \eqref{eq:odd_sod} and \eqref{eq:odd_sod2} using \eqref{eq:p_sod}:
$$D^b(P)=\left\langle \otheta^{-1} \mb A,\otheta^{-1}\mb B,\mb A,\mb B',\mb A(1),\mb B'(1),\otheta \mb A(1),\otheta \mb B(1)\right\rangle.$$
Here and below, we suppress pullbacks for brevity, and Serre twists are with respect to $\mc O_\pi(1)$. Tensoring the leftmost two blocks by $\omega_P^{-1}=\otheta(2)$ gives
$$D^b(P)=\left\langle\mb A,\mb B',\mb A(1),\mb B'(1),\otheta \mb A(1),\otheta \mb B(1),\mb A(2),\mb B(2)\right\rangle.$$
To continue, we need
\begin{lem}\label{lem:hecke_ortho}
We have
\begin{equation}\label{eq:hecke_ortho1}
    \otheta \mb A(1)\subset {}^\perp \mb A(2)
\end{equation}
\begin{equation}\label{eq:hecke_ortho2}
    \otheta\mb B(1)\subset {}^\perp \left\langle \mb A(2),\mb B(2)\right\rangle.
\end{equation}
\end{lem}
This allows us to reorder: 
$$D^b(P)=\left\langle\mb A,\mb B',\mb A(1),\mb B'(1),\mb A(2),\otheta \mb A(1),\mb B(2),\otheta \mb B(1)\right\rangle.$$
Finally, moving the rightmost block to the left proves
\begin{prp}\label{prp:hecke_sod_ab}
We have a semiorthogonal decomposition
$$D^b(P)=\left\langle\mb B(-1),\mb A,\mb B',\mb A(1),\mb B'(1),\mb A(2),\otheta \mb A(1),\mb B(2)\right\rangle.$$
\end{prp}
We now introduce counterparts of these blocks in $D^b(\eN)$. Let
\begin{align*}
    C_k&=\bLambda^{2-g+k} \boxpow{\eF}{g-2-2k}\\
    D_k&=\bLambda^{2-g+k} \boxpow{\eF}{g-1-2k}.
\end{align*}
\begin{lem}\label{lem:cd}
The Fourier--Mukai functors $\mc P_{C_k}$ (for $0\leq k\leq \floor{\frac{g-2}{2}}$) and $\mc P_{D_k}$ (for $0\leq k\leq \floor{\frac{g-1}{2}}$) into $D^b(\eN)$ are fully faithful. There are full subcategories $\mb C,\mb D,\mb D'\subset D^b(\eN)$ with semiorthogonal decompositions
\begin{align*}
    \mb C &= \langle C_k\rangle_{0\leq k\leq \floor{\frac{g-2}{2}}}\\
    \mb D &= \langle D_k\rangle_{1\leq k\leq \floor{\frac{g-1}{2}}}\\
    \mb D' &= \langle D_k\rangle_{0\leq k\leq \floor{\frac{g-1}{2}}},
\end{align*}
ordered by decreasing $k$.
\end{lem}
\begin{lem}\label{lem:cd_quasibps}
$\mb C\subset \mbb B_{2-g}$ and $\mb D,\mb D'\subset\mbb B_{3-g}$.
\end{lem}
In particular, by \cref{prp:quasibps}, the images of $\mb C,\mb D,\mb D'$ under $\sigma^*$ are admissible subcategories of $D^b(P)$, which we denote by the same symbols. We can now state the main result of this section.
\begin{thm} \label{thm:hecke_sod_cd}
We have semiorthogonal decompositions
\begin{align*}
    \mbb B_{2-g}&=\left\langle \etheta^{-1}\mb C,\mb C,\etheta \mb C,\etheta^2\mb C\right\rangle\\
    \mbb B_{3-g}&=\left\langle \etheta^{-1}\mb D,\mb D',\etheta\mb D',\etheta^2\mb D\right\rangle.
\end{align*}
 Moreover, we have a semiorthogonal decomposition $D^b(P)=\langle \sigma^*\mbb B_{2-g},\sigma^* \mbb B_{3-g}\rangle$. Its refinement in terms of $C_k$ and $D_k$ is related by a mutation to the refinement of \cref{prp:hecke_sod_ab} in terms of $A_k$ and $B_k$.
\end{thm}
\begin{rem}\label{rem:quasibps_sods}
Per \cref{rem:quasibps}, twisting by $\bLambda$ gives a semiorthogonal decomposition of $\mbb B_w$ for any $w\in \mbb Z$, and since $\omega_P^{-1}=\bLambda\etheta^2$, we have $D^b(P)=\langle\sigma^*\mbb B_w,\sigma^*\mbb B_{w+1}\rangle$. When $g$ is odd, our semiorthogonal decomposition of $\mbb B_1$ is of the form conjectured in \cite{belmans} (but note that $\mbb B_1$ is a \textit{twisted} noncommuative resolution).
\end{rem}
\begin{rem}\label{rem:hecke_bars}
    Recall the notation $\boxpow{\bar{\mc E}}{a}$ from \cref{ntn:stablepairs}. We write $\bar{A}_k=\otheta^{2-g+k}\boxpow{\bar{\mc E}}{g-2-2k}$, and similarly for $\bar{B}_k$, $\bar{C}_k$, $\bar{D}_k$. By the mutation of \cite{tevelevbraid}*{Theorem 6.14}, we can alternatively write $\mb A=\langle\bar A_k\rangle_{0\leq k\leq \floor{\frac{g-2}{2}}}$ ordered by increasing $k$, and similarly for $\mb B$ and $\mb B'$.
\end{rem}
\begin{lem}\label{lem:cd_bars}
There are mutations $\langle C_{\floor{\frac{g-2}{2}}},\dots,C_{k+1},C_k\rangle\to \langle \bar C_k,\bar C_{k+1},\dots,\bar C_{\floor{\frac{g-2}{2}}}\rangle$ for $0\leq k\leq \floor{\frac{g-2}{2}}$ and $\langle D_{\floor{\frac{g-1}{2}}},\dots,D_{k+1},D_k\rangle\to \langle \bar D_k,\bar D_{k+1},\dots,\bar D_{\floor{\frac{g-1}{2}}}\rangle$ for $0\leq k\leq \floor{\frac{g-1}{2}}$.
\end{lem}

\begin{figure}
    \centering
    \[\begin{tikzcd}[cramped,row sep=large,column sep=small]
    \bar B_{p-i}&\bar C_{p-i}(1)&\bar C_{p-i+1}(1)\dots\bar C_{p}(1)&D_q\dots D_{p-i+1}&A_{p-i-1}(1)&\\
    \bar C_{p-i}(1)&\bar D_{p-i}&&\,&&\\
    &&\,&\,&D_{p-i}&A_{p-i-1}(1)\\
    &&\,&\,&C_{p-i-1}(1)&D_{p-i}\\
    \bar C_{p-i-1}(1)&\bar C_{p-i}(1)&\bar C_{p-i+1}(1)\dots\bar C_{p}(1)&D_q\dots D_{p-i+1}&D_{p-i}&
    \arrow[from=1-5,to=3-6,no head,bend left=30]
    \arrow[from=1-1,to=2-2,no head,bend right=5]
    \arrow[from=1-2,to=2-1,no head,crossing over,bend right = 10]
    \arrow[from=3-6,to=4-5,no head,bend left=20]
    \arrow[from=3-5,to=4-6,no head,crossing over,bend left=20]
    \arrow[from=4-6,to=5-5,no head, bend left=20]
    \arrow[from=4-5,to=5-1,no head,bend right=10]
    \arrow[from=1-3,to=5-3,no head,crossing over]
    \arrow[from=2-1,to=5-2,no head,crossing over,bend right]
    \arrow[from=2-2,to=3-5,no head,bend right=10]
    \arrow[from=1-4,to=3-4,no head, crossing over]
    \arrow[from=2-4,to=5-4,no head]
    \end{tikzcd}\]
    \caption{The inductive step.}
    \label{fig:hecke_mutation}
\end{figure}
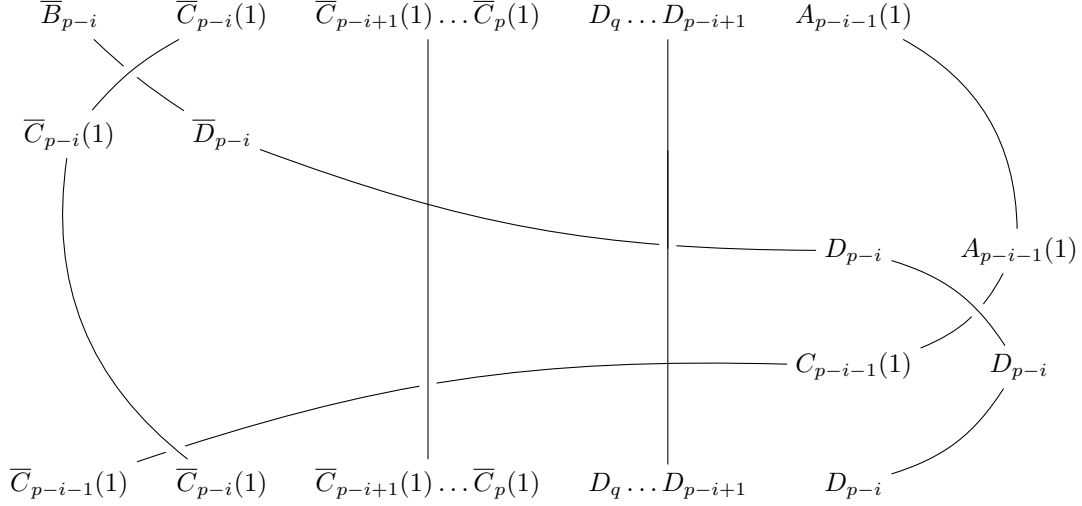

\begin{lem}\label{lem:basic_hecke}
Given that $\langle \bar{B}_\ell\rangle\subset \langle \bar{C}_\ell(1)\rangle^\perp$ with $0\leq \ell\leq \floor{\frac{g-2}{2}}$, or $\langle D_\ell\rangle\subset \langle A_{\ell-1}(1)\rangle^\perp$ with $1\leq \ell \leq \floor{\frac{g-1}{2}}$, there is a mutation
\begin{equation}\label{eq:basic_hecke1}
\begin{tikzcd}
    \bar{B}_\ell\arrow[dr,no head] & \bar{C}_\ell(1)\arrow[dl,no head,crossing over]\\
    \bar{C}_\ell(1) & \bar{D}_\ell
\end{tikzcd}
\end{equation}
or
\begin{equation}\label{eq:basic_hecke2}
\begin{tikzcd}
    D_\ell & A_{\ell-1}(1)\arrow[dl,no head]\\
    C_{\ell-1}(1) & D_\ell \arrow[ul,no head, crossing over]
\end{tikzcd}    
\end{equation}
respectively in $D^b(P)$.
\end{lem}
\begin{proof}[Proof of \cref{thm:hecke_sod_cd}]
We first mutate $\langle \mb B,\mb A(1)\rangle$ into $\langle \mb C(1),\mb D\rangle$. Let $p=\floor{\frac{g-2}{2}}$ and $q=\floor{\frac{g-1}{2}}$.
\begin{clm}
    For $0\leq i \leq p$,
    \begin{equation}\label{eq:hecke_ind}
        \langle \bar B_1,\dots,\bar B_{p-i}, \bar C_{p-i}(1),\dots,\bar C_{p}(1), D_q,\dots,D_{p-i+1}, A_{p-i-1}(1),\dots, A_0(1) \rangle 
    \end{equation}
    is a semiorthogonal decomposition of $\langle\mb B,\mb A(1)\rangle$. Each is a mutation of the previous.

\end{clm}
We prove the claim by induction on $i$. When $i=0$, there are two cases: If $g$ is even, so $q=p$, then \eqref{eq:hecke_ind} reads
$\langle \bar B_1,\dots,\bar B_p,C_p(1),A_{p-1}(1),\dots,A_0(1)\rangle$. Since $A_p=\otheta^{-p}=\bLambda^{-p}=\bar C_p$ by \eqref{eq:thetaislambda}, this is exactly $\langle \mb B,\mb A(1)\rangle$ after introducing bars as in \cref{rem:hecke_bars}. If $g$ is odd, so $q=p+1$, then \eqref{eq:hecke_ind} reads $\langle \bar B_1,\dots,\bar B_p,C_p(1),D_{p+1},A_{p-1}(1),\dots,A_0(1)\rangle$. Since $\bar B_{p+1}=D_{p+1}$, this is obtained from $\langle \mb B,\mb A(1)\rangle$ by the mutation \eqref{eq:basic_hecke2}.

For the inductive step, we refer to \cref{fig:hecke_mutation}. We first apply \eqref{eq:basic_hecke1} in the top left. By \cref{lem:cd_quasibps,prp:quasibps}, $\langle \bar D_{p-i}\rangle$ is contained in ${}^\perp\mb C(1)$, so we may move it past $\langle\bar C_{p-i+1}(1),\dots,\bar C_{p}(1)\rangle$ (note that $ \sigma^*\mbb B_w(1)=\sigma^*\mbb B_w$ by \eqref{eq:thetaisO1} and \cref{rem:quasibps}). Using \cref{lem:cd_bars}, we mutate $$\langle\bar D_{p-i},D_q,\dots,D_{p-i+1}\rangle\to\langle\bar D_{p-i},\bar D_{p-i+1},\dots,\bar D_{q}\rangle\to \langle D_q,\dots,D_{p-i+1},D_{p-i}\rangle.$$
The second half of the figure is similar, this time using \eqref{eq:basic_hecke2}. This proves the claim.

Taking $i=p$ gives the desired mutation of $\langle \mb B,\mb A(1)\rangle$ into $\langle \mb C(1),\mb D\rangle$. We can similarly mutate $\langle \mb B',\mb A(1)\rangle=\langle \bar B_0,\mb B,\mb A(1)\rangle\to\langle \bar B_0,\mb C(1),\mb D\rangle \to \langle \mb C(1),\mb D,D_0\rangle=\langle \mb C(1),\mb D'\rangle$ using the top half of \cref{fig:hecke_mutation} with $i=p$. Since $A_k^\dual=\otheta^{g-2} \bar A_k$, $B_k^\dual=\otheta^{g-3}\bar B_k$, $C_k^\dual=\bLambda^{g-2} \bar C_k$, and $D_k^\dual=\bLambda^{g-3} \bar D_k$ (cf. \eqref{eq:dualbar}), we have $\langle\otheta \mb A(1),\mb B(2)\rangle^\dual=\otheta^{g-3}\langle \mb B(-2),\mb A(-1)\rangle\to \bLambda^{g-3}\langle \mb C(-1),\mb D(-2)\rangle=\langle \mb D(2),\bLambda C(1)\rangle$. Hence there is a mutation of $\langle\otheta \mb A(1),\mb B(2)\rangle$ into $\langle \mb D(2),\bLambda\mb C(1)\rangle$ dual to the one described above. Put together, we can mutate \cref{prp:hecke_sod_ab} into 
$$D^b(P)=\left\langle\mb C,\mb D(-1),\mb C(1),\mb D',\mb C(2),\mb D'(1),\mb D(2),\bLambda \mb C(1)\right\rangle.$$
We tensor the last block by $\omega_P=\otheta^{-1}(-2)=\bLambda^{-1}(-2)$, moving it to the left:
$$D^b(P)=\left\langle\mb C(-1),\mb C,\mb D(-1),\mb C(1),\mb D',\mb C(2),\mb D'(1),\mb D(2)\right\rangle.$$
Finally, we reorder using \cref{prp:quasibps,lem:cd_quasibps} to obtain
$$D^b(P)=\left\langle\mb C(-1),\mb C,\mb C(1),\mb C(2),\mb D(-1),\mb D',\mb D'(1),\mb D(2)\right\rangle.$$ 
It follows that $$\sigma^*\mbb B_{2-g}=\langle\mb C(-1),\mb C,\mb C(1),\mb C(2)\rangle=\langle \etheta^{-1}\mb C,\mb C,\etheta \mb C,\etheta^2\mb C\rangle$$ and $$\sigma^*\mbb B_{3-g}=\langle\mb D(-1),\mb D',\mb D'(1),\mb D(2)\rangle=\langle \etheta^{-1}\mb D,\mb D',\etheta\mb D',\etheta^2\mb D\rangle,$$ as claimed.
\end{proof}

\begin{proof}[Proof of \cref{lem:hecke_ortho}]
By \cite{tevelevtorresbgmn}*{Lemma 10.2}, to prove \eqref{eq:hecke_ortho1} we must show that $$R\Hom_P\left(\otheta^{k+1}\boxpow{\oE}{g-2-2k}_D, \otheta^{\ell}\boxpow{\oE}{g-2-2\ell}_{D'}(1)\right)=0$$
for $0\leq k,\ell\leq \floor{\frac{g-2}{2}}$ and $D\in \Sym^{g-2-2k} C$, $D'\in \Sym^{g-2-2\ell} C$. Since $\pi_* \mc O_\pi(1)=\oE_q=\otheta\oE_q^\dual$, this is equivalent to 
$$R\Hom_{\oN}\left(\otheta^{k}\boxpow{\oE}{g-2-2k}_D, \otheta^{\ell}\boxpow{\oE}{g-2-2\ell}_{D'}\otimes \oE_q^\dual\right)=0$$
by the projection formula. As in \cites{tevelevtorresbgmn,tevelevbraid}, there is a birational morphism $\zeta:\hat M\to \oN$, where $\hat M=M_{g-1}(2g-1)$ from \cref{ntn:stablepairs}. We have $\zeta^*\boxpow{\oE}{g-2-2k}_D=Z^{g-2-2k}\boxpow{\oF}{g-2-2k}_D$ and $\zeta^*\otheta=Z^2 \hat\bLambda$, where $\oF$ is the universal bundle on $\hat M$, $Z=\mc O(1,g-1)$, and $\hat\bLambda=\mc O(0,-1)$. $\zeta^*$ is fully faithful, so it remains to show that
$$R\Gamma_{\hat M}\left(\hat\bLambda^{\ell-k}Z^{-1}(\boxpow{\oF}{g-2-2k}_D\otimes\oF_q)^\dual\otimes \boxpow{\oF}{g-2-2\ell}_{D'}\right)=0.$$
This follows from \cite{tevelevtorresbgmn}*{Theorem 4.1, Remark 4.2}. Indeed, $-1-2k<\ell-k<1+2\ell$. Likewise, \eqref{eq:hecke_ortho2} follows from
\begin{equation}\label{eq:hecke_ortho3}
    R\Gamma_{\hat M}\left(\hat\bLambda^{\ell-k}Z^{-2}(\boxpow{\oF}{g-1-2k}_D\otimes\oF_q)^\dual\otimes \boxpow{\oF}{g-2-2\ell}_{D'}\right)=0
\end{equation}
(for $1\leq k\leq \floor{\frac{g-1}{2}}$ and $0\leq \ell \leq\floor{\frac{g-2}{2}})$) and
\begin{equation}\label{eq:hecke_ortho4}
    R\Gamma_{\hat M}\left(\hat\bLambda^{\ell-k}Z^{-1}(\boxpow{\oF}{g-1-2k}_D\otimes\oF_q)^\dual\otimes \boxpow{\oF}{g-1-2\ell}_{D'}\right)=0
\end{equation}
(for $1\leq k,\ell\leq \floor{\frac{g-1}{2}}$). \eqref{eq:hecke_ortho4} follows immediately from \cite{tevelevtorresbgmn}*{Theorem 4.1}, since $-2k<\ell-k<2\ell$. After Serre duality with $\omega_{\hat M}=Z^{-3}\hat\bLambda^{-2}$, \eqref{eq:hecke_ortho3} becomes 
$$R\Gamma_{\hat M}\left(\hat\bLambda^{k-\ell-2}Z^{-1}(\boxpow{\oF}{g-2-2\ell}_{D'})^\dual\otimes\boxpow{\oF}{g-1-2k}_D\otimes\oF_q \right)=0.$$
 Since $-2-2\ell<k-\ell-2<2k-1$, we're done.
\end{proof}

\begin{proof}[Proof of \cref{lem:cd}]
Instead of $\eN$, we work on the isomorphic moduli stack $\eN'$ of semistable rank $2$ bundles of fixed determinant of degree $2g$ (this was called $\eN$ in \cref{sec:plainweave}). It is enough to prove the claim with $\eF$ replaced by the universal bundle on $C\times \eN'$, since they correspond up to an irrelevant twist. We use the notation of \eqref{eq:diagram}. By \cref{cor:windows}, it will suffice to check that pullback along $\alpha:\mc Z^\circ\to \eN'$ is fully faithful on $D^b(\eN')_w$. Since $\mc A$ and $\mc B$ have $\mbb{G}_m$-weight $1$, \cref{prp:Zcirckoszul} implies that the weight-$0$ component of $\alpha_*\mc O_{\mc Z^\circ}$ is $\mc O_{\eN'}$. Then for $X,Y\in D^b(\eN')_w$, we have $\Hom(\alpha^*X,\alpha^*Y)=\Hom(X,Y\otimes \alpha_*\mc O_{\mc Z^\circ})=\Hom(X,Y)$.
\end{proof}
\begin{proof}[Proof of \cref{lem:cd_quasibps}]
Objects in $\mb C$ and $\mb D,\mb D'$ have weights $2-g$ and $3-g$ respectively under scalar automorphisms, so we need only check \eqref{eq:quasibps} when $\nu$ corresponds to the one-parameter subgroup $\lambda$ stabilizing a split bundle considered in the proof of \cref{thm:main_precise}. As there, $n_\nu=g-1$, $\operatorname{weight}\nu^*\bLambda=-1$, and objects in $\nu^*\langle C_k\rangle$ and $\nu^*\langle D_k\rangle$ have weights in the ranges $[k,g-2-k]\subseteq [\frac{-1}{2},g-\frac{3}{2}]$ and $[k-1,g-2-k]\subseteq [-1,g-2]$, respectively (cf. \cref{prp:weights}). 
\end{proof}
\begin{proof}[Proof of \cref{lem:cd_bars}]
This follows from \cref{lem:makebars} and the proof of \cref{lem:cd}.
\end{proof}
\begin{proof}[Proof of \cref{lem:basic_hecke}]
After clearing powers of $\otheta=\bLambda$, \eqref{eq:basic_hecke1} reads
$$\begin{tikzcd}
    \boxpow{\bar\oE}{g-1-2\ell}\arrow[dr,no head] & \boxpow{\bar \eF}{g-2-2\ell}(1)\arrow[dl,no head,crossing over]\\
    \boxpow{\bar \eF}{g-2-2\ell}(1) & \boxpow{\bar\eF}{g-1-2\ell}
\end{tikzcd}$$
and \eqref{eq:basic_hecke2} is dual to
$$\begin{tikzcd}
    \boxpow{\bar\oE}{g-2\ell}\arrow[dr,no head] & \boxpow{\bar \eF}{g-1-2\ell}(1)\arrow[dl,no head,crossing over]\\
    \boxpow{\bar \eF}{g-1-2\ell}(1) & \boxpow{\bar\eF}{g-2\ell}
\end{tikzcd},$$
so it will suffice to give a mutation 
$$\begin{tikzcd}
    \boxpow{\bar\oE}{k}\arrow[dr,no head] & \boxpow{\bar \eF}{k-1}(1)\arrow[dl,no head,crossing over]\\
    \boxpow{\bar \eF}{k-1}(1) & \boxpow{\bar\eF}{k}
\end{tikzcd}$$
where $1\leq k\leq g-1$. Semiorthogonality at the top is assumed, and at the bottom it follows from \cref{prp:quasibps,lem:cd_quasibps}.
Therefore it remains to give an exact triangle 
\begin{equation}\label{eq:basic_hecke_triangle}
\boxpow{\bar\eF}{k}\to \boxpow{\bar\oE}{k}\to (i_{k-1})_*\boxpow{\bar\eF}{k-1}(1)\to
\end{equation}
in $D^b(\Sym^k C\times P)$, where $i_{k-1}:\Sym^{k-1} C\times P\hookrightarrow \Sym^k C\times P$ is the inclusion of the image of $q\times C^{k-1}\times P$ under the quotient map. When $k=1$, this is just \eqref{eq:heckeuniversal}. For $k=2$, recall that $\boxpow{\bar\eF}{2}$ is the equivariant pushfoward to $\Sym^2 C\times P$ of the $S_2$-equivariant sheaf $\pi_1^*\eF\otimes \pi_2^*\eF \otimes \mathrm{sgn}$ on $C^2\times P$, where $\pi_i:C^2\times P\to C\times P$ are the projections. Since $\eF$ is quasi-isomorphic to the complex $\oE\to \mc O_{q\times P}(1)$, $\pi_1^*\eF\otimes \pi_2^*\eF$ is equivariantly quasi-isomorphic to the complex
$$\pi_1^*\oE\otimes \pi_2^*\oE\to (\pi_1^*\oE|_{C\times q\times P}\oplus \pi_2^*\oE|_{q\times C\times P})(1)\to \mc O_{q\times q\times P}(2)\otimes \mathrm{sgn}$$
(note that $\mc O_{q\times C\times P}(1)\otimes^L \mc O_{C\times q\times P}(1)=\mc O_{q\times q\times P}(2)$ by transversality; the factor of $\mathrm{sgn}$ appears for degree reasons). Hence there is a morphism $\pi_1^*\eF\otimes \pi_2^*\eF\otimes\mathrm{sgn}\to \pi_1^*\oE\otimes \pi_2^*\oE\otimes\mathrm{sgn}$ with cone $$(\pi_1^*\oE|_{C\times q\times P}\oplus \pi_2^*\oE|_{q\times C\times P})(1)\otimes \mathrm{sgn} \to \mc O_{q\times q\times P}(2),$$ whose equivariant pushforward to $\Sym^2 C\times P$ is $[i_*\oE(1)\to i_*\mc O(2)]\simeq i_*\eF(1)$. This proves \eqref{eq:basic_hecke_triangle} for $k=2$.

Similarly, for any $k$, $\bigotimes_{i=1}^k\pi_i^*\eF$ is $S_k$-equivariantly quasi-isomorphic to the complex 
$$\bigotimes_{i=1}^k\pi_i^*[\oE\to \mc O_{q\times P}(1)]=[\oE_{k,k}\to \oE_{k,k-1}\to\cdots\to \oE_{k,0}]$$ with
$$\oE_{k,p}=\bigoplus_{|I|=p}(\bigotimes_{i\in I}\pi_i^*\oE)|_{C^I}(k-p),$$ 
where the direct sum is over subsets $I\subseteq \{1,\dots,k\}$, $C^I$ is the corresponding copy of $C^{p}\subseteq C^k$ (e.g., $C^{\{1,\dots,p\}}=C^{p}\times \{q\}^{k-p}$), and the $S_k$ action is induced from the $S_p\times S_{k-p}$ action on $(\bigotimes_{i=1}^p\pi_i^*\oE)|_{C^p\times \{q\}^{k-p}}(p)$ where $S_p$ acts by permuting the tensor factors and $S_{k-p}$ acts by its sign character (once again, transversality ensures that there are no contributions from higher $\mathrm{Tor}$s). By Frobenius reciprocity, taking $S_k$-invariants of $\mc E_{k,p}\otimes \mathrm{sgn}$ for $p<k$ is the same as taking $S_p\times S_{k-p}$-invariants of $(\bigotimes_{i=1}^p\pi_i^*\oE)|_{C^p\times \{q\}^{k-p}}(p)$, where now $S_{k-p}$ acts trivially. By restricting to $S_p\times S_{k-1-p}\times S_1$ and inducing to $S_{k-1}\times S_1$, we find that the $S_k$-equivariant pushfoward to $\Sym^k C\times P$ of $$\mathrm{Cone}(\bigotimes_{i=1}^k\pi_i^*\eF\otimes \mathrm{sgn}\to\bigotimes_{i=1}^k\pi_i^*\oE\otimes\mathrm{sgn})=[\oE_{k,k-1}\to\cdots\to \oE_{k,0}]\otimes \mathrm{sgn}$$ is quasi-isomorphic to the $S_{k-1}$-equivariant pushforward of $$[\oE_{k-1,k-1}\to\cdots\to\oE_{k-1,0}](1)\otimes \mathrm{sgn}\simeq(\bigotimes_{i=1}^{k-1}\pi_i^*\eF)(1)\otimes\mathrm{sgn}\in D^b_{S_{k-1}}(C^{k-1}\times P)$$ to $i_{k-1}(\Sym^{k-1}C\times P)$,
which is precisely $(i_{k-1})_*\boxpow{\bar\eF}{k-1}(1)$.
\end{proof}
\section{Basic weaving patterns}\label{sec:basicweaving}
In this section, we prove the various technical lemmas used in \cref{sec:stablepairssod}. The main statements and their proofs are taken from \cite{tevelevbraid} (where they are proved only for $d=2g-1$) with minor modifications and some additional details. We include them here mainly to verify that the numerical bookkeeping required to apply \cref{thm:hardvanishing} below remains valid in any degree $d\leq 2g$; we omit those proofs with no dependence on~$d$. Let $v=\floor{\frac{d-1}{2}}$ and $M_i=M_i(d)$. We write $\mc P_K$ for the Fourier--Mukai functor with kernel $K$.

\subsection{Cross Warp}\label{sub:crosswarp}
Closely following \cite{tevelevbraid}*{Section 4}, we prove the following:
\begin{thm}[Basic Cross Warp, cf. \cite{tevelevbraid}*{Theorem 3.2}]\label{thm:crosswarp} For $0\leq k\leq i\leq v$, we have:
\begin{thmlist}
    \item \label{thm:F_ff} $\cP_{\boxpow{\cF^\dual}{k}}:D^b(\Sym^k C)\to D^b(M_i)$ is fully faithful.
    \item \label{thm:D_ff} $\cP_{\cD^k_i}:D^b(\Sym^k C)\to D^b(M_i)$ is fully faithful.
    \item \label{thm:D_windows} If $k\leq i-1$, then $\iota\langle\cD^k_{i-1}\rangle=\langle\cD^k_i\rangle$ where $\iota$ is the windows embedding of \cref{prp:windows}. Moreover, objects in $\langle\cD^k_{i-1}\rangle$ descend from objects with weights in the range $[0,k]$ for this wall crossing. 
    \item \label{thm:crosswarp_mutation}There is an admissible subcategory of $D^b(M_i)$ with semiorthogonal decompositions 
    $$\langle\boxpow{\cF^\dual}{k-1},\dots,\cO,\cD^k_i,\cD^{k-1}_i\bLambda^{-1},\dots,\bLambda^{-k}\rangle$$ and $$\langle\cD^{k-1}_i\bLambda^{-1},\dots,\bLambda^{-k},\boxpow{\cF^\dual}{k},\boxpow{\cF^\dual}{k-1},\dots,\cO\rangle$$
    related by the mutation in \cref{fig:crosswarp}.

\end{thmlist}
\end{thm}
\begin{rem} \label{rem:F_windows}
    It follows from \cite{tevelevtorresbgmn}*{Section 3} that $\iota\langle \boxpow{\cF^\dual}{k}\rangle=\langle \boxpow{\cF^\dual}{k}\rangle\subset D^b(M_i)$ for $k<i$, and that objects in this subcategory have weights in the range $[0,k]$.
\end{rem}
\begin{cor}\label{cor:Fs_contain_D}
    $\langle \cD^k_i\rangle$ is contained in the subcategory generated by $\langle \boxpow{\cF^\dual}{\ell}\bLambda^{-m}\rangle$ with $0\leq \ell\leq k$, $0\leq m\leq k-\ell$.
\end{cor}
\begin{proof}
    \cref{thm:crosswarp_mutation} and induction on $k$.
\end{proof}
\begin{proof}[Proof of \cref{thm:F_ff}]
    Note that since $\boxpow{\cF^\dual}{k}\cong (\Lambda^{-k}\boxtimes\bLambda^{-k})\otimes \boxpow{\cF}{k}$ (see \cite{tevelevbraid}*{Lemma 3.4}), we may equivalently prove that $\cP_{\boxpow{\cF}{k}}$ is fully faithful. The case $d\leq 2g-1$ is \cite{tevelevtorresbgmn}*{Theorem 9.2}, and the proof for $d=2g$ is the same. Indeed, since \cite{tevelevtorresbgmn}*{Cor. 8.2, Cor. 8.4, Thm. 9.6} are proved for all $d\leq 2g+1$, we only need to verify that $k\leq 3g-2i-2$ for $k\leq i$. Since $k\leq g-1$, we have $3g-2i-2\geq g$, so this holds.
\end{proof}

The other parts of \cref{thm:crosswarp} will be proved shortly. We require the vanishing theorems for tensor vector bundles proved in \cite{tevelevtorresbgmn}.
\begin{thm}[\cite{tevelevbraid}*{Theorems 3.8, 3.9, 3.10}]\label{thm:hardvanishing}
Let $d'$, $j$, $a$, $b$, and $t$ be integers with $2<d'\leq 2g+1$, $1\leq j\leq \floor{\frac{d'-1}{2}}$, and let $D\in\Sym^a C$, $D'\in\Sym^bC$.
\begin{thmlist}[itemsep=1em]
    \item\label{thm:hv1} 
        If $a,\ b\leq d'+g-2j-1$, $t\notin [0,a]$, and $a-j-1 < t < d'+g-2j-1-b,$ then
        \begin{equation*}
            R\Gamma\left(M_j(d'),(\boxpow{\bar\cF}{a}_D)^\dual\otimes \boxpow{\bar\cF}{b}_{D'}\otimes\bLambda^t\right)=0
        \end{equation*}
        and, if $D=\sum\alpha_i x_i$,
        \begin{equation*}
            R\Gamma\left(M_j(d'),\bigotimes_i{(\cF^\dual_{x_i}})^{\otimes \alpha_i}\otimes\boxpow{\bar\cF}{b}_{D'}\otimes\bLambda^t\right)=0.
        \end{equation*}

    \item\label{thm:hv2} 
        If $a < t < d'+g-2j-1-b$, then 
        \begin{equation*}
        R\Gamma\left(M_j(d'),(\boxpow{\bar\cF}{a}_D)^\dual\otimes \boxpow{\bar\cF}{b}_{D'}\otimes\bLambda^t\right)=
        R\Gamma\left(M_j(d'),(\boxpow{\cF}{a}_D)^\dual\otimes\boxpow{\cF}{b}_{D'} \otimes\bLambda^t\right)=
        0.
        \end{equation*}
        Moreover, the same vanishing holds with $j=0$ for any $d'>0$.
        
    \item\label{thm:hv3} 
        If $a\le j$, $b< d'+g-2j-1$, and $D\not\le D'$ (for example, $a>b$), then 
        \begin{equation*}
        R\Gamma\left(M_j(d'),(\boxpow{\bar\cF}{a}_D)^\dual\otimes \boxpow{\bar\cF}{b}_{D'}\right)=0.
        \end{equation*}
\end{thmlist}
\end{thm}

\begin{rem}
    \cref{thm:hv3} is not completely proved in \cites{tevelevtorresbgmn,tevelevbraid} for $d'=2g,2g+1$. It comes from \cite{tevelevtorresbgmn}*{Theorem 9.6}, which is proved under the inductive assumption $\cP_{\boxpow{\cF}{k}}$ is fully faithful for $k<a$. This is true for $d'=2g$ by \cref{thm:F_ff}, and for $d'=2g+1$ by the same argument (since then $a\leq i\leq g$, so $k\leq g-1=d'+g-2g-2$).
\end{rem}

\begin{cnj}\label{cnj:hardervanishing}
    \cref{thm:hardvanishing} holds for any $d'>2$ and $j\leq v$ with $3j\leq d'+g-1$. 
\end{cnj}
\begin{rem}
We expect this conjecture would follow from a careful analysis of the proofs given in \cite{tevelevtorresbgmn}. The main purpose of the hypothesis $d'\leq 2g+1$ is to ensure that $3j\leq d'+g-1$ (needed for calculations using windows) for all $j\leq \floor{(d'-1)/2}$; one needs to check that it is not needed otherwise (e.g., for inequalities related to ample cones). 
\end{rem}

\begin{lem}[cf. \cite{tevelevbraid}*{Lemma 4.1}] \label{lem:subcats}
    For $1\leq k\leq i$, $D^b(M_i)$ contains admissible subcategories
    $$\langle \boxpow{\cF^\dual}{k-1},\dots,\cF^\dual,\cO,\bLambda^{-k}\rangle \text{ and } \langle\bLambda^{-k},\boxpow{\cF^\dual}{k},\dots,\cF^\dual,\cO\rangle.$$
\end{lem}
\begin{proof}
    All blocks are admissible by \cref{thm:F_ff}. By \cite{tevelevtorresbgmn}*{Lemma 10.2}, it suffices to check semi-orthogonality on skycraper sheaves of closed points $D\in \Sym^a C$, $D'\in \Sym^b C$. For $b<a\leq k$, we have $b<v< d+g-2i-1$ since $v< \frac{d+g-1}{3}$. Thus
    $$R\Hom(\boxpow{\cF^\dual}{b}_{D'},\boxpow{\cF^\dual}{a}_{D})=R\Hom(\boxpow{\bar\cF}{a}_{D},\boxpow{\bar\cF}{b}_{D'})=0$$
    by \cite{tevelevbraid}*{Lemma 3.4} and \cref{thm:hv3}. Next, for $0\leq b\leq k$, we have $$0, l\leq i<d+g-2i-1,\qquad -k\not\in [0,0],\qquad -i-1<-k<d+g-2i-1-b,$$
    so
    $$R\Hom(\boxpow{\cF^\dual}{b}_{D'},\bLambda^{-k})=R\Gamma(M_i,\boxpow{\bar\cF}{b}_{D'}\otimes\bLambda^{-k})=0$$
    by \cref{thm:hv1}. Finally, for $0\leq a<k$, we have $a<k\leq v<d+g-2i-1$, so 
    $$R\Hom(\bLambda^{-k},\boxpow{\cF^\dual}{a}_D)=R\Gamma(M_i,(\boxpow{\bar\cF}{a}_D)^\dual\otimes\bLambda^k)=0$$
    by \cref{thm:hv2}.
\end{proof}

We introduce the following complexes on $\Sym^k C\times M_i$, which will allow us to break the mutation of \cref{thm:crosswarp_mutation} into two stages (\cref{lem:FtoFbullet,lem:FbullettoD}):

\begin{ntn}[\cite{tevelevbraid}*{Definitions 4.2 and 4.5}]\label{ntn:Fbullet}
    Let $\cF^\bullet=[\cF^\dual\to\cO]\in D^b(C\times M_i)$, with $\cO$ in degree $0$ and the map given by contraction with the universal section of $\cF$. Let $$\boxpow{\cF^\bullet}{k}=\tau_*^{S_k}\left(\pi_1^*\cF^\bullet{\otimes}\cdots{\otimes}\pi_k^*\cF^\bullet\otimes \mathrm{sgn}\right)$$
    where $\tau: C^k\to \Sym^k C$ is the quotient, $\pi_\ell:C^k\to C$ are the projections, and $\mathrm{sgn}$ is the sign representation of $S_k$.
\end{ntn}

We prove parts (b), (c), and (d) \cref{thm:crosswarp} by induction on $k$. The following lemmas are proved for each $k$ as part of the same induction.

\begin{lem}\label{lem:cw_orthogonal}
    For $0\leq \ell\leq k$, $0\leq m<k$ and $X$, $Y$ skyscraper sheaves at $D\in \Sym^\ell C$, $D'\in \Sym^m C$, respectively,
    \begin{equation}\label{eq:cw_orthogonal1}
        R\Hom(\cP_{\boxpow{\cF^\dual}{\ell}}(X),\cP_{\cD^m_i\bLambda^{m-k}}(Y))=0
    \end{equation}
    and
    \begin{equation}\label{eq:cw_orthogonal2}
        R\Hom(\cP_{\cD^\ell_i\bLambda^{\ell-k}}(X),\cP_{\boxpow{\cF^\dual}{m}}(Y))=0.
    \end{equation}
\end{lem}
\begin{proof}
    We have 
    $$R\Gamma\left(M_{i-m}(\Lambda(-2D')),\left(\boxpow{\cF^\dual}{\ell}_D\right)^\dual\bLambda^{m-k}\right)=R\Gamma\left(M_{i-m}(d-2m),\boxpow{\bar\cF}{\ell}_D\bLambda^{m-k}\right)=0$$
    by \cref{thm:hv1}. Indeed, we have $$0,\ell\leq v< d+g-2i-1,\qquad m-k\not\in [0,0],\qquad -i+m-1<m-k<d+g-2i-1-\ell.$$
    This proves \eqref{eq:cw_orthogonal1}. For \eqref{eq:cw_orthogonal2}, we have $\omega_{M_i}|_{M_{i-\ell}(d-2\ell)}=\cO(-3,-d-g+4+3\ell)=\omega_{M_{i-\ell}(d-2\ell)}\bLambda^{-\ell}$, so
    \begin{align*}
        R\Hom\left(\boxpow{\cF^\dual}{m}_{D'},\cO_{M_{i-\ell}(d-2\ell)}\bLambda^{\ell-k}\omega_{M_i}\right)
        &=R\Gamma\left(M_{i-\ell}(d-2\ell),\boxpow{\bar\cF}{m}_{D'}\bLambda^{-k}\omega_{M_{i-\ell}(d-2\ell)}\right)\\
        &=R\Gamma\left(M_{i-\ell}(d-2\ell),\left(\boxpow{\bar\cF}{m}_{D'}\right)^\dual\bLambda^{k}\right)[\dots].
    \end{align*}
    by Serre duality (we suppress the shift in degree). This vanishes by \cref{thm:hv2}, as $m<k<3g-2i-1$ (note that either $j>0$ and $d'>2$, or $j=0$ and $d'>0$). \eqref{eq:cw_orthogonal2} then follows from Serre duality.
\end{proof}

\begin{rem}\label{rem:cw_orthogonal}
    Since $\cP_{\boxpow{\cF^\dual}{\ell}}$ and $\cP_{\cD^m_i\bLambda^{m-k}}(Y))$ are fully faithful by \cref{thm:F_ff} and the inductive hypothesis, this proves \eqref{eq:cw_orthogonal1} for any $X\in D^b(\Sym^\ell C)$, $Y\in D^b(\Sym^\ell C)$ by \cite{tevelevtorresbgmn}*{Lemma 10.2}. The same is true for \eqref{eq:cw_orthogonal2} with $\ell<k$.
\end{rem}

\begin{lem}[cf. \cite{tevelevbraid}*{Corollary 4.13}]\label{lem:cw_triangle}
    There exist $G_k,H_k\in D^b(\Sym^k C\times M_i)$ such that for all $X\in D^b(\Sym^k C)$, we have $\cP_{G_k}(X)\in \langle\boxpow{\cF^\dual}{k-1},\dots,\cF^\dual,\cO\rangle$, $\cP_{H_k}(X)\in\langle \cD^{k-1}_i\bLambda^{-1},\dots,\cD^1_i\bLambda^{1-k},\bLambda^{-k}\rangle$, and exact triangles
    \begin{equation}\label{eq:cw_triangle1}
        \cP_{G_k}(X)\to \cP_{\boxpow{\cF^\bullet}{k}}(X)\to \cP_{\boxpow{\cF^\dual}{k}}(X)[k]\to
    \end{equation}
    and
    \begin{equation}\label{eq:cw_triangle2}
        \cP_{H_k}(X)\to \cP_{\boxpow{\cF^\bullet}{k}}(X)\to \cP_{\cD^k_i}(X(-B/2))\to 
    \end{equation}
    where $B\subset\Sym^k C$ is the diagonal divisor.
\end{lem}
\begin{proof}
    The proofs of Lemma 4.7 through Corollary 4.13 in \cite{tevelevbraid} carry over to any $d$ without change. 
\end{proof}

\begin{lem}[cf. \cite{tevelevbraid}*{Lemma 4.14}] \label{lem:FtoFbullet}
    $\cP_{\boxpow{\cF^\bullet}{k}}$ is fully faithful and there is a mutation
    $$\langle\boxpow{\cF^\dual}{k},\boxpow{\cF^\dual}{k-1},\dots,\cO\rangle\to \langle\boxpow{\cF^\dual}{k-1},\dots,\cO,\boxpow{\cF^\bullet}{k}\rangle.$$
\end{lem}
\begin{proof}
    Applying $R\Hom(-,\cP_{\boxpow{\cF^\dual}{m}}(Y))$ to \eqref{eq:cw_triangle2} and using \eqref{eq:cw_orthogonal2} gives \begin{equation}\label{eq:FtoFbullet1}
    R\Hom(\cP_{\boxpow{\cF^\bullet}{k}}(X),\cP_{\boxpow{\cF^\dual}{m}}(Y))=0
    \end{equation} 
    for $m<k$ and $X,Y$ skyscraper sheaves (note $R\Hom(\cP_{H_k}(X),\cP_{\boxpow{\cF^\dual}{m}}(Y))=0$ by \cref{rem:cw_orthogonal}). Then
    \begin{align*}
        R\Hom(\cP_{\boxpow{\cF^\bullet}{k}}(X),\cP_{\boxpow{\cF^\bullet}{k}}(Y)) &= R\Hom(\cP_{\boxpow{\cF^\bullet}{k}}(X),\cP_{\boxpow{\cF^\dual}{k}}(Y)[k])\\
        &= R\Hom(\cP_{\boxpow{\cF^\dual}{k}}(X)[k],\cP_{\boxpow{\cF^\dual}{k}}(Y)[k])=R\Hom(X,Y)
    \end{align*}
    by \eqref{eq:cw_triangle1}, \eqref{eq:FtoFbullet1}, and \cref{lem:subcats}. By the Bondal--Orlov criterion, this proves full faithfulness. \cref{lem:subcats} and \eqref{eq:FtoFbullet1} give semiorthogonality, and \eqref{eq:cw_triangle1} yields the claimed mutation.
\end{proof}
\begin{lem}\label{lem:FbullettoD}
    $\cP_{\cD^k_i}$ is fully faithful and there is a mutation $$\langle\cD^{k-1}_i\bLambda^{-1},\dots,\bLambda^{-k},\boxpow{\cF^\bullet}{k}\rangle\to\langle\cD^k_i,\cD^{k-1}_i\bLambda^{-1},\dots,\bLambda^{-k}\rangle.$$
\end{lem}
\begin{proof}
    We show
    \begin{equation}\label{eq:FbullettoD1}
        R\Hom(\cP_{\boxpow{\cF^\bullet}{k}}(X),\cP_{\cD^m_i\bLambda^{m-k}}(Y))=0
    \end{equation}
    and 
    \begin{equation}\label{eq:FbullettoD2}
        R\Hom(\cP_{\cD^m_i\bLambda^{m-k}}(Y),\cP_{\cD^k_i}(X))=0
    \end{equation}
    for $m<k$ and $X,Y$ skyscraper sheaves. \eqref{eq:FbullettoD1} follows from \eqref{eq:cw_orthogonal1} and \eqref{eq:cw_triangle1}. By the inductive hypothesis \cref{thm:crosswarp_mutation} and induction on $m$, it suffices to prove \eqref{eq:FbullettoD2} that 
    $$0=R\Hom\left(\boxpow{\cF^\dual}{b}_{D'}\bLambda^{-t},\cO_{M_{i-k}(d-2k)}\right)=R\Gamma\left(M_{i-k}(d-2k),\boxpow{\bar\cF}{b}_{D'}\bLambda^{t}\right)$$
    for $0\leq b<k$, $0<t\leq k-b$, $D\in \Sym^k C$, $D'\in \Sym^b C$ (cf. \cref{cor:Fs_contain_D}). Since $k\leq i<d+g-2i-1$, we have $0<t<d+g-2i-1-b$, so the required vanishing follows from \cref{thm:hv2}. 
    
    To prove full faithfulness, we compute
    \begin{align*}
        R\Hom(\cP_{\cD^k_i}(X),\cP_{\cD^k_i}(Y))
        &=R\Hom(\cP_{\boxpow{\cF^\bullet}{k}}(X(B/2)),\cP_{\cD^k_i}(Y))\\
        &=R\Hom(\cP_{\boxpow{\cF^\bullet}{k}}(X(B/2)),\cP_{\boxpow{\cF^\bullet}{k}}(Y(B/2)))=R\Hom(X,Y)
    \end{align*}
    by \eqref{eq:cw_triangle2}, \eqref{eq:FbullettoD2}, \eqref{eq:FbullettoD1}, and \cref{lem:FtoFbullet}. Semiorthonality then follows from \eqref{eq:FbullettoD1} and \eqref{eq:FbullettoD2}, together with the inductive hypothesis. Finally, \eqref{eq:cw_triangle2} gives the required mutation.
\end{proof}
\begin{proof}[Proof of \cref{thm:crosswarp}]
    We proved \ref{thm:F_ff} above and \ref{thm:D_ff} in \cref{lem:FbullettoD}. \ref{thm:crosswarp_mutation} follows from \cref{lem:cw_orthogonal,lem:FtoFbullet,lem:FbullettoD}, as depicted in \cref{fig:crosswarp_steps}. Finally, \ref{thm:D_windows} follows from \ref{thm:crosswarp_mutation} and the inductive hypothesis. Indeed, we have $\iota\langle \boxpow{\cF^\dual}{\ell}\rangle=\langle \boxpow{\cF^\dual}{\ell}\rangle$ for $0\leq\ell\leq k$ by \cref{rem:F_windows} and $\iota\langle\cD^m_{i-1}\bLambda^{m-k}\rangle=\langle\cD^m_{i}\bLambda^{m-k}\rangle$ by the inductive hypothesis (recall that $\bLambda$ has weight $-1$). Performing the mutation of \ref{thm:crosswarp_mutation}, embedding via $\iota$, and undoing the mutation on the other side shows that $\iota\langle\cD^k_{i-1}\rangle=\langle\cD^k_{i}\rangle$. It is clear from the bottom decomposition that the entire subcategory has weights in the range $[0,k]$, so this holds for $\langle\cD^k_{i-1}\rangle$ a fortiori.
\end{proof}
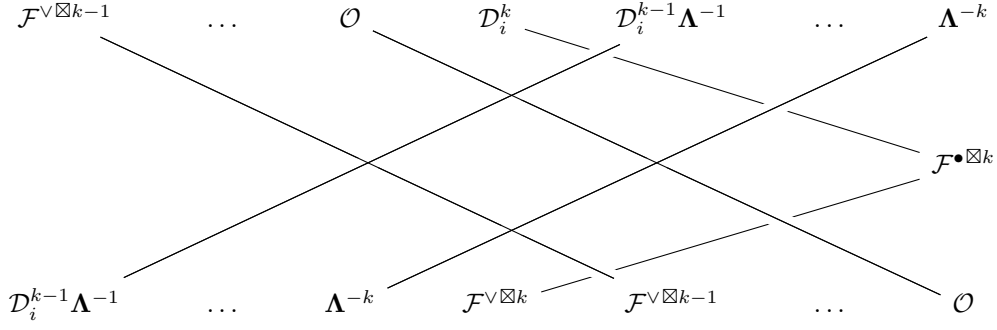
\begin{figure}
    \centering
    \[\begin{tikzcd}
    {\boxpow{\cF^\dual}{k-1}} & \dots & {\mathcal{O}} & {\mathcal{D}^k_i} & {\mathcal{D}^{k-1}_i\mathbf{\Lambda}^{-1}} & \dots & {\mathbf{\Lambda}^{-k}} \\
	\\ & & & & & & \boxpow{\cF^\bullet}{k}
    \\ 
	\\
	{\mathcal{D}^{k-1}_i\mathbf{\Lambda}^{-1}} & \dots & {\mathbf{\Lambda}^{-k}} & \boxpow{\cF^\dual}{k} & \boxpow{\cF^\dual}{k-1} & \dots & {\mathcal{O}}
    \arrow[from=1-4,to=3-7,no head]
    \arrow[from=3-7,to=5-4,no head]
	\arrow[no head, crossing over, from=1-3, to=5-7]
	\arrow[no head, crossing over, from=1-1, to=5-5]
	\arrow[no head, crossing over, from=1-5, to=5-1]
	\arrow[no head, crossing over, from=1-7, to=5-3]
 	\arrow[no head, from=1-3, to=5-7]
	\arrow[no head, from=1-1, to=5-5]
	\arrow[no head, from=1-5, to=5-1]
	\arrow[no head, from=1-7, to=5-3]
\end{tikzcd}\]
    \caption{The mutation of \cref{thm:crosswarp_mutation} in two steps, cf. \cite{tevelevbraid}*{Figure 9}.}
    \label{fig:crosswarp_steps}
\end{figure}


\subsection{Farey Twill}\label{sub:farey_twill}
We now turn to the orthogonalities and mutations required to implement the Farey Twill. Recall that for integers $k,s$ with $0\leq k\leq v$ and $t\in [k,v+1)$, we define $\cD^{k,s}_t=\cD^k_{\floor{t}}\otimes L^{k,s}_t$ where 
$$L^{k,s}_t=
\begin{cases}
    \cO(s,sk)                                                          &  k = \floor{t} \\
    \cO\left(\floor{\frac{s}{t-k}},s+\floor{\frac{s}{t-k}}(k-1)\right) &  k < \floor{t}. 
\end{cases}$$
\begin{lem}[cf. \cite{tevelevbraid}*{Lemma 2.7}]\label{lem:ft_windows}
    For $k\leq i-1$ and $\epsilon\ll 0$, we have $\iota\langle\cD^{k,s}_{i-\epsilon}\rangle=\langle\cD^{k,s}_{i}\rangle$.
\end{lem}
\begin{proof}
    When $i=1$, $\iota$ is simply pullback under $M_1\to M_0$, so $\iota\langle\cD^{0,s}_{1-\epsilon})\rangle=\langle\iota(\cO(s))\rangle=\langle\cD^{0,s}_1\rangle$. Suppose $i>1$. By \cite{tevelevtorresbgmn}*{Remark 3.22}, $\cO(m,n)$ has weight $n+m(1-i)$ for the wall crossing $M_{i-1}\dashrightarrow M_i$. Hence $L^{k,s}_{i-\epsilon}=L^{k,s}_{i}$ has weight $s-\floor{\frac{s}{i-k}}(i-k)\in [0,i-k)$. Hence the objects in $\langle(\cD^{k,s}_{i-\epsilon})\rangle$ have weights in the range $[0,i)$. The lemma then follows from \cref{prp:windows,thm:D_windows}.
\end{proof}
\begin{lem}[cf. \cite{tevelevbraid}*{Lemma 2.8}]\label{lem:ft_orthogonal}
    Let $x=\frac{s}{t-k}\not\in \ZZ$. If $k<k'$ and either $x=\frac{s'}{t-k'}$ or $s'=t-k'=0$,  then
    \begin{equation}\label{eq:ft_orthogonal1}
        R\Hom(\cP_{\cD^{k,s}_t}(X),\cP_{\cD^{k',s'}_t}(Y))=0
    \end{equation}
    for any $X\in D^b(\Sym^k C)$, $Y\in D^b(\Sym^{k'} C)$.
\end{lem}
\begin{proof}
    By \cref{cor:Fs_contain_D}, it suffices to prove
    $$R\Hom(\cP_{\boxpow{\cF^\dual}{\ell}\bLambda^{-m}L^{k,s}_t}(Z),\cP_{\cD^{k',s'}_t}(Y))=0$$
    for $0\leq\ell\leq k$, $0\leq m\leq k-\ell$, and $Z\in D^b(\Sym^\ell C)$. Without loss of generality, we take $Y$ and $Z$ to be skyscraper sheaves, reducing \eqref{eq:ft_orthogonal1} to
    
    \begin{equation}\label{eq:ft_orthogonal2}
        R\Hom\left(\boxpow{\cF^\dual}{\ell}_D \bLambda^{-m}L^{k,s}_t,\cO_{M_0(\Lambda(-2D'))}L^{k',s'}_t\right)=0.
    \end{equation}
    
    We first consider the case $k'=\floor{t}$, so $L^{k',s'}_t=\cO(s',s'k')$. By \cite{tevelevtorresbgmn}*{Remark 3.7}, $\cO(m,n)$ restricts to $\cO(n+m(1-k'))=\bLambda^{-n-m(1-k')}$ on the fibers $M_0(d-2k')\cong\PP^{d+g-2-2k'}$ of $D^{k'}_{k'}$ over $\Sym^{k'} C$ . Noting that $s-s'=x(k'-k)$, we see that $(L^{k,s}_t)^{-1}L^{k',s'}_t$ restricts to $\bLambda^{\{x\}(k'-k)}$ where $\{x\}=x-\floor{x}$. The case $k'<\floor{t}$ is similar: $\cO(m,n)$ restricts to $\cO(m,n-mk')$ on $M_{i-k'}(d-2k')$, so $(L^{k,s}_t)^{-1}L^{k',s'}_t$ restricts to $\cO(0,\{x\}(k-k'))=\bLambda^{\{x\}(k'-k)}$. Either way, \eqref{eq:ft_orthogonal2} becomes
    $$R\Gamma\left(M_{\floor{t}-k'}(d-2k'),\boxpow{\bar\cF}{\ell}_D \bLambda^{m+\{x\}(k'-k)}\right)=0.$$
    This follows from \cref{thm:hv2} once we verify 
    $$0<m+\{x\}(k'-k)<d-2\floor{t}+g-1-\ell.$$
    Since $\{x\}>0$, the first inequality is clear. For the second, we check $m+\ell+\{x\}(k'-k)+2\floor{t}<d+g-1$. Since $$m+\ell\leq k,\quad \{x\}(k'-k)<k'-k,\quad k'\leq \floor{t},\text{ and}\quad 3\floor{t}\leq 3v<d+g-1,$$ we're done. 
\end{proof}
\begin{lem}[cf. \cite{tevelevbraid}*{Lemma 2.9}]\label{lem:ft_mutation}
    For all $k\leq i$, there is a mutation
    $$\langle \cD^k_i,\dots,\cD^{i-1}_i,\cD^i_i\rangle\to \langle\cD^{i}_i,\cD^{i-1}_i(-1,2-i),\dots,\cD^k_i(-1,1-k)\rangle.$$
\end{lem}
\begin{proof}
    Semiorthogonality on the left-hand side follows from \cref{prp:windows,thm:D_windows}. As in \cite{tevelevbraid}*{Claim 3.9}, we have for every $X\in \Sym^n C$ an exact triangle
    \begin{equation} \label{eq:ft_mutation1}
        \cP_{\cD^n_i}(X)(-1,1-k)\to\cP_{\cD^n_i}(X)\to\cP_{\cO_{E^n_i}}(X)\to
    \end{equation}
    where $\cP_{\cO_{E^n_i}}(X)\in \langle \cD^{n+1}_i,\dots,\cD^i_i\rangle$. Here, $E^n_i$ is the divisor $\{(D,F,s)\in D^n_i:\deg Z(s)\geq n+1\}$ where $Z(s)$ denotes the scheme of zeros. We proceed by downward induction on $k$, the case $i=k$ being trivial. By the inductive hypothesis, it suffices to give a mutation $\langle \cD^k_i,\cD^{k+1}_i,\dots,\cD^i_i\rangle\to \langle \cD^{k+1}_i,\dots,\cD^i_i,\cD^k_i(-1,1-k)\rangle$. By \eqref{eq:ft_mutation1}, we are left to show that $\cD^k_i(-1,1-k)$ is contained in ${}^\perp\langle \cD^{k+1}_i,\dots,\cD^i_i\rangle$. Again by the inductive hypothesis, this amounts to showing that $$R\Hom(\cP_{\cD^k_i}(X)(-1,1-k),\cP_{\cD^n_i}(Y)(-1,1-n))=0$$
    for $k<n\leq i$, $X\in D^b(\Sym^k C)$, $Y\in \Sym^n C$. (Note that $\cO(-1,1-i)$ has trivial restriction to $M_0(d-2i)$, so $\cD^i_i(-1,1-i)\cong \cD^i_i$.) By \cref{cor:Fs_contain_D}, the required vanishing is
    $$R\Gamma(M_{i-n}(d-2n),\boxpow{\bar\cF}{\ell}_D\bLambda^{m+n-k})=0$$
    for $0\leq \ell\leq k$, $0\leq m\leq k-\ell$, $D\in \Sym^\ell C$. We have $0<m+n-k<d+g-2i-1-\ell$, so we're done by \cref{thm:hv2}.
\end{proof}

\begin{lem}\label{lem:ft_reordering}
    For $k,k',j,j'\geq 0$ with $j'+k'<j+k\leq i$, we have
    \begin{equation*}
        R\Hom(\cP_{\cD^k_i\bLambda^{-j}}(X),\cP_{\cD^{k'}_i\bLambda^{-j'}}(Y))=0
    \end{equation*}
    for any $X\in D^b(\Sym^k C)$, $Y\in D^b(\Sym^{k'} C)$.
\end{lem}
\begin{proof}
    As above, it suffices to show that
    \begin{equation*}\label{eq:ft_reordering1}
        R\Hom(\cO_{M_{i-k}(d-2k)}\bLambda^{-j},\boxpow{\cF^\dual}{\ell'}_{D}\bLambda^{-j'-m'})=0
    \end{equation*}
    for $0\leq\ell'\leq k'$, $0\leq m'\leq k'-\ell'$, $D\in \Sym^{\ell'}C$. By Serre duality applied twice (cf. the proof of \cref{lem:cw_orthogonal}), this is equivalent to 
    \begin{equation*}
        R\Gamma\left(M_{i-k}(d-2k),(\boxpow{\bar\cF}{\ell'}_D)^\dual\bLambda^{j-j'-m'+k}\right)=0.
    \end{equation*}
    This follows from \cref{thm:hv2} once $\ell'<j-j'-m'+k<d+g-2i-1$. For the first, we have $(j-j'+k)-(m'+\ell')>k'-k'=0$; for the second, we have $j-j'-m'+k\leq j+k\leq i<d+g-2i-1$.
\end{proof}

\subsection{Broken Loom} Finally, we adapt the reordering trick from \cite{tevelevbraid}*{Section 5}. We recall the following lemma:
\begin{lem}[\cite{tevelevbraid}*{Lemma 5.4, Remark 5.7}]\label{lem:bl_vanishing}
    Suppose $2<d'\leq 2g+1$ and $1\leq j\leq \floor{\frac{d-1}{2}}$. Let $D\in \Sym^a C$, $D'\in \Sym^b C$ with $a,b\leq j$ (\cite{tevelevbraid} has $a,b\leq \min(j,d'+g-2j-1)$, but $j\leq d'+g-2j-1$ already), and let $t$ be an integer with $a-j-1<t<d'+g-2j-1-b$ and $2t<a-b$. Then
    \begin{equation}\label{eq:bl_vanishing1}
        R\Gamma\left(M_j(d'),\left(\boxpow{\bar\cF}{a}_{D}\right)^\dual \otimes \boxpow{\bar\cF}{b}_{D'}\otimes \bLambda^t\right)=0.
    \end{equation}
\end{lem}

\begin{lem}\label{lem:bl_reordering}
    Let $\lambda,\lambda',k,k'$ be integers with $k,k',\lambda-2k,\lambda'-2k'\geq 0$ and $\lambda-k,\lambda'-k'\leq i_d$. If $\lambda<\lambda'$, then
    $$R\Hom\left(\bLambda^{-k}\boxpow{\cF^\dual}{\lambda-2k}_D,\bLambda^{-k'}\boxpow{\cF^\dual}{\lambda'-2k'}_{D'}\right)=0$$
    for any $D\in \Sym^{\lambda-2k} C$, $D'\in \Sym^{\lambda'-2k'} C$.
\end{lem}
\begin{proof}
    We need $$R\Gamma\left(M_{i_d}(d),(\boxpow{\bar\cF}{\lambda-2k}_{D})^\dual \boxpow{\bar\cF}{\lambda'-2k'}_{D'}\bLambda^{k-k'}\right)=0,$$ which is exactly \eqref{eq:bl_vanishing1}. We clearly have $\lambda-2k,\lambda'-2k'\leq i_d$, and $2(k-k')<\lambda'-2k'-(\lambda-2k')$ since $\lambda'>\lambda$. It remains to show that $$\lambda'-2k'-i_d-1<k-k'<d+g-2i_d-1-\lambda+2k.$$ We have $\lambda'-k'-i_d-1\leq-1< k$ for the first inequality, and $d+g-2i_d-1-(\lambda-k)\geq d+g-3i_d-1>0\geq -k$ for the second.
\end{proof}

\bibliography{main}

\end{document}